\documentclass{article}

\usepackage{arxiv}

\usepackage[utf8]{inputenc} % allow utf-8 input
\usepackage[T1]{fontenc}    % use 8-bit T1 fonts
\usepackage{hyperref}       % hyperlinks
\usepackage{url}            % simple URL typesetting
\usepackage{booktabs}       % professional-quality tables
\usepackage{amsfonts}       % blackboard math symbols
\usepackage{nicefrac}       % compact symbols for 1/2, etc.
\usepackage{microtype}      % microtypography
\usepackage{lipsum}		% Can be removed after putting your text content
\usepackage{graphicx}
\usepackage{natbib}
\usepackage{doi}
\usepackage{enumerate}
\usepackage{enumitem}
\usepackage{amsmath, amsthm, amssymb, algorithm, algpseudocode}
\usepackage[table]{xcolor}
\usepackage{tikz}
\DeclareMathOperator{\HC}{HC}  % Harish-Chandra category
\DeclareMathOperator{\Ext}{Ext}
\DeclareMathOperator{\Hom}{Hom}

\DeclareMathOperator{\wt}{wt}  % weight
\newcommand{\HCfin}{\HC_{\text{fin}}}
\newcommand{\FF}{\mathbb{F}}
\newcommand{\KK}{\mathbb{K}}
\newcommand{\ZZ}{\mathbb{Z}}  % Added missing command
\newcommand{\QQ}{\mathbb{Q}}  % Added missing command

\newcommand{\cA}{\mathcal{A}}

\newcommand{\cT}{\mathcal{T}}
\newcommand{\CC}{\mathbb{C}}
\newcommand{\Aut}{\operatorname{Aut}}

\newcommand{\Weyl}[1]{A\!\left(\FF[#1]\right)}

\newcommand{\Witt}[1]{W\!\left(\FF[#1]\right)}

\newcommand{\ev}{\operatorname{ev}}
\newcommand{\GL}{\operatorname{GL}}
\newcommand{\rank}{\operatorname{rank}}
\newtheorem{theorem}{Theorem}[section]

\newtheorem{definition}[theorem]{Definition}
\newtheorem{example}[theorem]{Example}

\newtheorem{problem}[theorem]{Problem}
\newcommand\mystyle{\everymath{\displaystyle}}
\mystyle
\title{Weyl-Type and Witt-Type Algebras with Exponential Generators:Structure, Automorphisms, and Representation Theory}

\author{\href{https://orcid.org/0000-0002-3816-5287}{\includegraphics[scale=0.06]{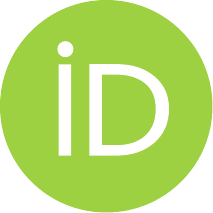}\hspace{1mm}M.H.M.~Rashid}\thanks{Corresponding Author} \\
	Department of Mathematics\&Statistics\\Faculty of Science P.O.Box(7)\\
	Mutah University University\\
	Mutah-Jordan \\
	\texttt{mrash@mutah.edu.jo}
}

\renewcommand{\shorttitle}{\textit{arXiv} Template}

\hypersetup{
pdftitle={Weyl-Type and Witt-Type Algebras with Exponential Generators:Structure, Automorphisms, and Representation Theory},
pdfsubject={q-bio.NC, q-bio.QM},
pdfauthor={M.H.M.Rashid},
pdfkeywords={Weyl-type algebras, Witt-type algebras, exponential generators, additive subgroups, automorphism groups, quantum deformations, weight modules, Harish-Chandra modules, Galois descent, category $\mathcal{O}$}}

\begin{document}
\maketitle

\begin{abstract}
	 This paper introduces and systematically studies a new class of non-commutative algebras—Weyl-type and Witt-type algebras—generated by differential operators with exponential and generalized power function coefficients. We define the expolynomial ring $R_{p,t,\cA} = \FF\bigl[ e^{\pm x^{p} e^{t}},\; e^{\cA x},\; x^{\cA} \bigr]$ associated to an additive subgroup $\cA \subset \FF$, and investigate its Ore extension $A_{p,t,\cA} = R_{p,t,\cA}[\partial; \delta]$ (Weyl-type) and its derivation algebra $\mathfrak{g}_{p,t,\cA} = \operatorname{Der}_{\FF}(R_{p,t,\cA})$ (Witt-type). Our main results establish: (1) the automorphism group of $R_{p,t,\cA}$ is isomorphic to $(\FF^{\times})^{2r+1} \rtimes \GL(2r+1,\ZZ)$; (2) a Galois descent theorem showing that fixed-point subalgebras under finite Galois actions recover the original Weyl-type algebra; (3) the non-existence of finite-dimensional simple modules for $A_{p,t,\cA}$; (4) the Zariski density of isomorphism classes in moduli spaces as transcendental parameters vary; (5) the stability of simplicity under generic quantum deformation; and (6) a complete representation-theoretic framework including the classification of irreducible weight modules, the construction of Harish–Chandra modules with BGG-type resolutions, and the structure of category $\mathcal{O}$. These results unify and extend classical theories of Weyl algebras, Witt algebras, and generalized Weyl algebras, while opening new directions in deformation theory, non-commutative geometry, and the representation theory of infinite-dimensional algebras.
\end{abstract}

% keywords can be removed
\keywords{Weyl-type algebras\and Witt-type algebras\and  exponential generators\and  additive subgroups\and automorphism groups\and quantum deformations\and weight modules\and Harish-Chandra modules\and  Galois descent\and category $\mathcal{O}$}
%%%%%%%%%%%%%%%%%%%%%%%%%%%%%%%%%%%%%%%%%%%%%%%%%%%%%%%%%%%%%%%%%%%%%%%%%%%%%%%%%%%%%%%%
%================================================================================
%=====================================================================================
\section{Introduction}
%=====================================================================================

The study of non-commutative algebras arising from differential operators and their symmetries has been a central theme in algebra, representation theory, and mathematical physics for over a century. Classical examples include the Weyl algebras, enveloping algebras of Lie algebras, and algebras of differential operators on algebraic varieties. More recently, algebras generated by differential operators with exponential or generalized power coefficients have attracted considerable attention due to their rich structure and connections with integrable systems, quantum groups, and deformation theory.

\subsection{Scholarly Contributions}
This paper builds upon several important lines of research in non-commutative algebra. The foundational work of Dixmier \cite{Dixmier68} on Weyl algebras established fundamental results about their simplicity and representation theory. Bavula \cite{Bavula92} introduced generalized Weyl algebras as a unifying framework for many algebras of differential operators. The study of Witt-type algebras with exponential coefficients was pioneered by Kawamoto \cite{Kawamoto86} and has connections with earlier work by Ree \cite{Ree56}. The representation theory follows the philosophy of weight modules and category $\mathcal{O}$ as developed for infinite-dimensional Lie algebras by Kac \cite{Kac74}, Rudakov \cite{Rudakov69}, and others. Our work extends these classical constructions by systematically incorporating exponential and power generators parameterized by additive subgroups, and by studying their Galois descent, automorphism groups, quantum deformations, and representation theory in a unified manner. For more details, the author refer to \cite{Angulo2022}.

\subsection{Applications}
The algebras studied in this paper have natural interpretations in several areas of mathematics and mathematical physics:
\begin{itemize}
    \item \textbf{Mathematical Physics:} They arise as algebras of symmetry operators for certain linear differential equations with exponential coefficients, and as deformation quantizations of Poisson algebras on algebraic tori.
    \item \textbf{Representation Theory:} They provide natural examples of infinite-dimensional algebras with well-behaved categories of weight modules, generalizing the representation theory of semisimple Lie algebras.
    \item \textbf{Deformation Theory:} The quantum deformations studied here connect to $q$-analysis and quantum groups, while the transcendental parameters lead to interesting families in moduli spaces.
    \item \textbf{Algebraic Geometry:} The commutative subalgebras $R_{p,t,\cA}$ are coordinate rings of algebraic tori, and their automorphism groups relate to birational geometry.
\end{itemize}

\subsection{Core of Study}
The central objects of this paper are the \emph{expolynomial rings}
\[
R_{p,t,\cA} = \FF\bigl[ e^{\pm x^{p} e^{t}},\; e^{\cA x},\; x^{\cA} \bigr]
\]
and their associated non-commutative algebras:
\begin{itemize}
    \item The \emph{Weyl-type algebras} $A_{p,t,\cA} = R_{p,t,\cA}[\partial; \delta]$, which are Ore extensions by a derivation $\partial$.
    \item The \emph{Witt-type algebras} $\mathfrak{g}_{p,t,\cA} = \operatorname{Der}_{\FF}(R_{p,t,\cA})$, which are Lie algebras of derivations.
\end{itemize}
Here $\cA$ is a finitely generated additive subgroup of a characteristic zero field $\FF$, $p \in \cA$, and $t \in \FF$. These algebras interpolate between classical Weyl/Witt algebras (when $\cA=\ZZ$ and $p=0$) and more exotic algebras with transcendental parameters.

\subsection{Significance of the Study}
This work provides a systematic investigation of algebras generated by differential operators with exponential and power function coefficients. The main contributions include:
\begin{enumerate}
    \item \textbf{Structural Results:} We determine the automorphism groups of expolynomial rings (Theorem~\ref{thm:aut-group}), establish Galois descent principles (Theorem~\ref{thm15}), and prove simplicity results for quantum deformations (Theorem~\ref{thm19}).
    \item \textbf{Representation Theory:} We classify irreducible weight modules for Weyl-type algebras (Theorem~\ref{thm:irreducible-weight}), study Harish-Chandra modules for Witt-type algebras (Theorem~\ref{thm:HC-modules}), and establish BGG-type resolutions in category $\mathcal{O}$ (Theorem~\ref{thm:BGG-resolution}).
    \item \textbf{Geometric Insights:} We show that families of these algebras are Zariski-dense in appropriate moduli spaces (Theorem~\ref{thm18}), revealing the richness of the deformation theory.
    \item \textbf{Methodological Advances:} We develop techniques that combine commutative algebra (Laurent polynomial rings), Galois theory, deformation theory, and highest weight category methods in a novel context.
\end{enumerate}

The paper is organized as follows: Section 1 collects preliminaries and definitions. Section 2 studies Galois descent, automorphisms, and quantum deformations. Section 3 establishes structure theorems for expolynomial rings. Section 4 develops the representation theory of Weyl-type and Witt-type algebras. Section 5 presents open problems for further research.
%=====================================================================================
\section{Preliminaries}
%=====================================================================================

This section collects the basic definitions, notations, and known results that will be used throughout the paper. We work over a field $\FF$ of characteristic zero, usually assumed to be algebraically closed (e.g., $\FF = \CC$).

\subsection{Basic Algebraic Structures}

\begin{definition}[Additive Subgroup $\cA$]
A \emph{rank-$r$ additive subgroup} of $\FF$ is a finitely generated $\ZZ$-submodule $\cA \subset (\FF, +)$ of rank $r$ over $\ZZ$. That is, $\cA \cong \ZZ^r \oplus T$ where $T$ is a finite torsion group. In most examples, $\cA$ will be a free $\ZZ$-module of rank $r$ (e.g., $\ZZ[\sqrt{d}]$ for a square-free integer $d$).
\end{definition}

\begin{definition}[Generalized Exponential and Power Functions]
For an additive subgroup $\cA \subset \FF$, we denote by
\[
\FF[e^{\cA x}, x^{\cA}]
\]
the commutative $\FF$-algebra generated by formal symbols
\[
\{ e^{\alpha x} \mid \alpha \in \cA \} \quad \text{and} \quad \{ x^{\alpha} \mid \alpha \in \cA \}
\]
subject to the relations:
\begin{align*}
e^{\alpha x} e^{\beta x} &= e^{(\alpha+\beta)x}, \\
x^{\alpha} x^{\beta} &= x^{\alpha+\beta}, \\
e^{\alpha x} x^{\beta} &= x^{\beta} e^{\alpha x},
\end{align*}
for all $\alpha, \beta \in \cA$. We interpret $x^{\alpha}$ as a formal power function and $e^{\alpha x}$ as a formal exponential symbol.
\end{definition}

\begin{definition}[Expolynomial Ring]
For a chosen element $p \in \cA$ and a parameter $t \in \FF$, we enlarge the above ring by adjoining a central invertible generator $e^{\pm x^{p} e^{t}}$ (or more generally $e^{\pm x^{p} e^{t x}}$). The resulting commutative $\FF$-algebra is denoted
\[
R_{p,t,\cA} = \FF\bigl[ e^{\pm x^{p} e^{t}},\; e^{\cA x},\; x^{\cA} \bigr].
\]
We call $R_{p,t,\cA}$ an \emph{expolynomial ring}. It is a Laurent polynomial ring in several variables (see Theorem~\ref{thm:aut-group}).
\end{definition}

\subsection{Weyl-Type and Witt-Type Algebras}

\begin{definition}[Weyl-Type Algebra]
The \emph{Weyl-type algebra} associated to the data $(p,t,\cA)$ is the Ore extension
\[
A_{p,t,\cA} = \Weyl{e^{\pm x^{p} e^{t}},\; e^{\cA x},\; x^{\cA}}
          = R_{p,t,\cA}[\partial; \delta],
\]
where $\partial$ is a derivation acting on $R_{p,t,\cA}$ by $\partial(x^{\alpha}) = \alpha x^{\alpha-1}$ (formally) and $\partial(e^{\alpha x}) = \alpha e^{\alpha x}$. The algebra is presented by generators
\[
\widehat{x},\; \partial,\; \widehat{e^{\alpha x}}\;(\alpha\in\cA),\; \widehat{e^{x^{p} e^{t x}}}
\]
with relations coming from the commutation rules
\[
[\partial, \widehat{x}] = 1,\quad
[\partial, \widehat{e^{\alpha x}}] = \alpha \widehat{e^{\alpha x}},\quad
[\partial, \widehat{e^{x^{p} e^{t x}}}] = \widehat{e^{x^{p} e^{t x}}}\cdot \bigl(p x^{p-1} e^{t x} + t x^{p} e^{t x}\bigr),
\]
and all other commutators among the generators (except the one involving $\partial$ and $\widehat{x}$) being zero. When $\cA=\ZZ$ and $p=0$, this reduces to the usual first Weyl algebra $A_1(\FF)$.
\end{definition}

\begin{definition}[Witt-Type Algebra]
The \emph{Witt-type algebra} is the Lie algebra of derivations of $R_{p,t,\cA}$:
\[
\mathfrak{g}_{p,t,\cA} = \Witt{e^{\pm x^{p} e^{t}},\; e^{\cA x},\; x^{\cA}}
                      = \operatorname{Der}_{\FF}\bigl(R_{p,t,\cA}\bigr).
\]
It has a natural basis $\{ e^{\alpha x} x^{\beta} \partial \mid \alpha,\beta\in\cA \}$ with Lie bracket
\[
[f\partial, g\partial] = \bigl(f\partial(g) - g\partial(f)\bigr)\partial.
\]
When $\cA=\ZZ$ and we omit the exponential generators, we recover the classical Witt algebra $\mathsf{W}_1$.
\end{definition}

\subsection{Representation-Theoretic Notions}

\begin{definition}[Weight Modules]
Let $A = A_{p,t,\cA}$. An $A$-module $M$ is a \emph{weight module} if it decomposes as
\[
M = \bigoplus_{\lambda \in \FF} M_{\lambda},
\qquad
M_{\lambda} = \{ m \in M \mid \widehat{x}\cdot m = \lambda m \}.
\]
The set $\wt(M) = \{ \lambda \in \FF \mid M_{\lambda} \neq 0 \}$ is called the \emph{set of weights} of $M$.
\end{definition}

\begin{definition}[Harish-Chandra Modules]
For the Witt-type algebra $\mathfrak{g} = \mathfrak{g}_{p,t,\cA}$, a $\mathfrak{g}$-module $M$ is called a \emph{Harish-Chandra module} if
\begin{enumerate}[label=(\roman*)]
    \item $M$ is a weight module with respect to the Cartan subalgebra $\FF\partial$,
    \item each weight space $M_{\lambda}$ is finite-dimensional, and
    \item the set $\wt(M)$ lies in a finite union of $\cA$-cosets.
\end{enumerate}
The category of Harish-Chandra modules for $\mathfrak{g}$ is denoted $\HC(\mathfrak{g})$.
\end{definition}

\begin{definition}[Category $\mathcal{O}$]
For $A = A_{p,t,\cA}$, choose a decomposition $\cA = \cA^{+} \cup \{0\} \cup \cA^{-}$ (a \emph{triangular decomposition}). The category $\mathcal{O}$ consists of all $A$-modules $M$ such that
\begin{enumerate}[label=(\roman*)]
    \item $M$ is a weight module,
    \item $M$ is finitely generated, and
    \item for each $m\in M$, $\dim \FF[\widehat{e^{\alpha x}} \mid \alpha\in\cA^{+}]\cdot m < \infty$.
\end{enumerate}
\end{definition}

\subsection{Galois Theory and Deformations}

\begin{definition}[Galois Fixed-Point Subalgebra]
Let $\KK/\FF$ be a finite Galois extension with Galois group $G$. If $A$ is an algebra over $\KK$, the \emph{fixed-point subalgebra} is
\[
A^{G} = \{ a \in A \mid \sigma(a)=a \ \forall \sigma\in G \}.
\]
Classical Galois descent says that when $A = B \otimes_{\FF} \KK$, we have $A^{G} \cong B$.
\end{definition}

\begin{definition}[Quantum Deformation]
A \emph{quantum deformation} of an algebra $A$ over $\FF$ is a family of algebras $\{A_{q}\}$ over $\FF(q)$ or $\FF[[q-1]]$ such that $A_{1} \cong A$ and the multiplication depends analytically on $q$. In our context, the deformation is given by replacing the canonical commutation relation $[\partial,\widehat{x}]=1$ with the $q$-deformed relation
\[
\partial \widehat{x} - q \widehat{x} \partial = 1,
\]
while keeping all other relations unchanged.
\end{definition}

\subsection{Key Known Results}

We shall frequently use the following classical facts:

\begin{itemize}
    \item The Weyl algebra $A_n(\FF)$ has no finite-dimensional simple modules when $\operatorname{char}\FF = 0$ (Dixmier~\cite{Dixmier68}).
    \item The automorphism group of a Laurent polynomial ring $\FF[z_1^{\pm1},\dots,z_m^{\pm1}]$ is isomorphic to $(\FF^{\times})^{m} \rtimes \GL(m,\ZZ)$ (see~\cite{Hungerford74}).
    \item For a finite Galois extension $\KK/\FF$, Galois descent gives an equivalence between $\FF$-algebras and $\KK$-algebras with a compatible $G$-action (see~\cite{AtiyahMacdonald69}).
    \item The category $\mathcal{O}$ for semisimple Lie algebras is a highest weight category with BGG resolutions (see~\cite{Bjork79}).
    \item Generalized Weyl algebras provide a uniform framework for many algebras of differential operators (see Bavula~\cite{Bavula92}).
\end{itemize}

The algebras studied in this paper can be viewed as generalizations of both the Weyl algebra and the algebras of generalized exponentials considered in~\cite{Bavula92, Kirillov76}. The Witt-type algebras extend the classical Witt algebra $\mathsf{W}_1$ by allowing exponential coefficients, related to the algebras studied by Kawamoto~\cite{Kawamoto86} and Ree~\cite{Ree56}. The representation theory follows the philosophy of weight modules and category $\mathcal{O}$ as developed for infinite-dimensional Lie algebras (see~\cite{Kac74, Rudakov69}).

All these structures will be investigated in detail in the subsequent sections.
%%%%%%%%%%%%%%%%%%%%%%%%%%%%%%%%%%%%%%%%%%%%%%%%%%%%%%%%%%%%%%%%%%%%%
%%%%%%%%%%%%%%%%%%%%%%%%%%%%%%%%%%%%%%%%%%%%%%%%%%%%%%%%%%%%%%%%%%%%%
\section{Galois Descent, Automorphisms, and Quantum Deformations of Exponential Weyl Algebras}
%121212121212121212121212121212%121212121212121212121212121212
\begin{definition}[Galois Fixed-point Subalgebra]
    Let $A$ be an algebra over a Galois extension $\mathbb{K}/\mathbb{F}$ with Galois group $G$. The \emph{fixed-point subalgebra} under $G$ is
    \[
    A^G = \{ a \in A \mid \sigma(a) = a \ \forall \sigma \in G \}.
    \]
\end{definition}
%131313%131313%131313%131313%131313%131313%131313%131313
\begin{theorem}[Galois Descent for Algebraic Extensions]\label{thm15}
Let $\KK/\FF$ be a finite Galois extension with Galois group $G$. Then
\[
A\!\left(\KK[e^{\pm x^p e^{t}},\; e^{\cA x},\; x^{\cA}]\right)^G \cong \Weyl{e^{\pm x^p e^{t}},\; e^{\cA x},\; x^{\cA}},
\]
where the fixed-point algebra is taken with respect to the natural $G$-action induced from $\KK/\FF$.
\end{theorem}
\begin{proof}
Let $\KK/\FF$ be a finite Galois extension with Galois group $G$. Consider the Weyl-type algebra over $\KK$,
\[
A_\KK = A\!\left(\KK[e^{\pm x^p e^{t}},\; e^{\cA x},\; x^{\cA}]\right).
\]
The group $G$ acts on $\KK$ by field automorphisms, and this action extends naturally to an action on the polynomial-exponential ring
$\KK[e^{\pm x^p e^{t}},\; e^{\cA x},\; x^{\cA}]$ by acting on the coefficients from $\KK$ while leaving the variables $x_i$ and the exponential symbols
$e^{\pm x^p e^{t}}$, $e^{\cA x}$, $x^{\cA}$ invariant. More precisely, for $\sigma \in G$ and a monomial $\lambda f(x) \in \KK[e^{\pm x^p e^{t}},\; e^{\cA x},\; x^{\cA}]$
with $\lambda \in \KK$ and $f(x)$ a monomial in the $x_i$ and exponential factors, we define $\sigma(\lambda f(x)) = \sigma(\lambda) f(x)$.
This gives a $G$-action on the commutative ring by ring automorphisms.

The action further extends to the Weyl-type algebra $A_\KK$ by letting $G$ act on the coefficients of elements expressed in the $\KK$-basis
while fixing the generators $x_i$, $\partial_i$, $e^{\pm x^p e^{t}}$, $e^{\cA x}$, $x^{\cA}$. That is, if an element $a \in A_\KK$ is written as
a finite sum $\sum_{\alpha,\beta} c_{\alpha\beta} x^\alpha \partial^\beta$ with $c_{\alpha\beta} \in \KK[e^{\pm x^p e^{t}},\; e^{\cA x},\; x^{\cA}]$,
then $\sigma(a) = \sum_{\alpha,\beta} \sigma(c_{\alpha\beta}) x^\alpha \partial^\beta$. Since the commutation relations in $A_\KK$ involve only the variables
and the exponential symbols, which are fixed by $G$, and the action on coefficients is by ring automorphisms, each $\sigma$ is an automorphism of the $\KK$-algebra $A_\KK$.
Thus we have a $G$-action on $A_\KK$ by algebra automorphisms.

The fixed-point subalgebra is defined as
\[
A_\KK^G = \{ a \in A_\KK \mid \sigma(a) = a \text{ for all } \sigma \in G \}.
\]
We aim to show that $A_\KK^G$ is isomorphic to the Weyl-type algebra over $\FF$, namely $A_\FF = \Weyl{e^{\pm x^p e^{t}},\; e^{\cA x},\; x^{\cA}}$.

Observe that $A_\FF$ is naturally contained in $A_\KK$ via the inclusion $\FF \hookrightarrow \KK$. Since the $G$-action on $A_\KK$ fixes all elements
whose coefficients lie in $\FF$, we have $A_\FF \subseteq A_\KK^G$. To prove the reverse inclusion, take an arbitrary element $a \in A_\KK^G$.
Write $a$ as a finite sum
\[
a = \sum_{\alpha,\beta} f_{\alpha\beta}(x) x^\alpha \partial^\beta,
\]
where $f_{\alpha\beta}(x) \in \KK[e^{\pm x^p e^{t}},\; e^{\cA x},\; x^{\cA}]$. The condition $\sigma(a) = a$ for all $\sigma \in G$ implies that
for each pair $(\alpha,\beta)$, we have $\sigma(f_{\alpha\beta}(x)) = f_{\alpha\beta}(x)$. Hence each coefficient $f_{\alpha\beta}(x)$ is invariant under the $G$-action
on the ring $\KK[e^{\pm x^p e^{t}},\; e^{\cA x},\; x^{\cA}]$.

Now, by Galois descent for vector spaces, the fixed points of $G$ acting on $\KK[e^{\pm x^p e^{t}},\; e^{\cA x},\; x^{\cA}]$ are precisely
$\FF[e^{\pm x^p e^{t}},\; e^{\cA x},\; x^{\cA}]$. Indeed, consider the $\FF$-subspace $\FF[e^{\pm x^p e^{t}},\; e^{\cA x},\; x^{\cA}]$ of the $\KK$-vector space
$\KK[e^{\pm x^p e^{t}},\; e^{\cA x},\; x^{\cA}]$. Since $\KK/\FF$ is Galois, we have $\KK^G = \FF$, and the extension of scalars
$\KK \otimes_\FF \FF[e^{\pm x^p e^{t}},\; e^{\cA x},\; x^{\cA}] \cong \KK[e^{\pm x^p e^{t}},\; e^{\cA x},\; x^{\cA}]$ respects the $G$-action.
The fixed points of $G$ on the tensor product are exactly the elements of the form $1 \otimes g$ with $g \in \FF[e^{\pm x^p e^{t}},\; e^{\cA x},\; x^{\cA}]$.
Thus $f_{\alpha\beta}(x) \in \FF[e^{\pm x^p e^{t}},\; e^{\cA x},\; x^{\cA}]$ for all $\alpha,\beta$.

Consequently, $a$ can be expressed as a finite sum with coefficients in $\FF[e^{\pm x^p e^{t}},\; e^{\cA x},\; x^{\cA}]$, which means $a \in A_\FF$.
Therefore $A_\KK^G \subseteq A_\FF$. Combining the two inclusions yields $A_\KK^G = A_\FF$.

Hence we have an isomorphism of $\FF$-algebras
\[
A\!\left(\KK[e^{\pm x^p e^{t}},\; e^{\cA x},\; x^{\cA}]\right)^G \cong \Weyl{e^{\pm x^p e^{t}},\; e^{\cA x},\; x^{\cA}},
\]
where the left-hand side is the fixed-point subalgebra under the natural $G$-action induced from the Galois extension $\KK/\FF$.
\end{proof}
%EEEEEEEEEEEEEEEEEEEEEEEEEEEEEEEEEEEEEEEEEEEEEEEEEEEEEEEEEEEEEEEEEEEEEE
\begin{example}[Illustration of Theorem~\ref{thm15}]
Let $\FF = \mathbb{Q}$ and $\KK = \mathbb{Q}(\sqrt{2})$, which is a Galois extension of degree $2$ with Galois group $G = \{\mathrm{id}, \sigma\}$, where $\sigma(\sqrt{2}) = -\sqrt{2}$. Let $\cA = \mathbb{Z}[\sqrt{3}] = \{m + n\sqrt{3} \mid m, n \in \mathbb{Z}\}$, viewed as an additive subgroup of $\mathbb{C}$. Fix $p = 3$ and $t = \pi$, a transcendental number.

Consider the base algebra over $\KK$:
\[
R_\KK = \KK[\, e^{\pm x^3 e^{\pi x}},\; e^{\cA x},\; x^{\cA} ].
\]
Its elements are finite $\KK$-linear combinations of monomials
\[
e^{a x^3 e^{\pi x}} e^{\alpha x} x^{\beta}, \quad a \in \mathbb{Z},\; \alpha \in \cA,\; \beta \in \cA.
\]
The Galois group $G$ acts on $\KK$ by field automorphisms, and this action extends naturally to $R_\KK$ by acting on coefficients:
\[
\sigma\left( \sum_{\text{finite}} c_\gamma \, e^{a_\gamma x^3 e^{\pi x}} e^{\alpha_\gamma x} x^{\beta_\gamma} \right)
= \sum_{\text{finite}} \sigma(c_\gamma) \, e^{a_\gamma x^3 e^{\pi x}} e^{\alpha_\gamma x} x^{\beta_\gamma},
\]
where $c_\gamma \in \KK$.

Now form the Weyl-type algebra over $\KK$:
\[
A_\KK = A\!\left( R_\KK \right) = \KK\langle \widehat{e^{x^3 e^{\pi x}}},\; \widehat{e^{x}},\; \widehat{x},\; \partial \rangle / \mathcal{R},
\]
with the relations $\mathcal{R}$ given by the commutation rules
\[
[\partial, \widehat{x}] = 1, \quad [\partial, \widehat{e^{x}}] = \widehat{e^{x}}, \quad [\partial, \widehat{e^{x^3 e^{\pi x}}}] = \widehat{e^{x^3 e^{\pi x}}} \cdot \big( 3x^2 e^{\pi x} + \pi x^3 e^{\pi x} \big).
\]
The $G$-action on $A_\KK$ is defined coefficient-wise: for $P = \sum c_\gamma \, M_\gamma \in A_\KK$, where $M_\gamma$ are monomials in the generators,
\[
\sigma(P) = \sum \sigma(c_\gamma) \, M_\gamma.
\]
This action is well-defined because $\sigma$ preserves the relations $\mathcal{R}$ (since the coefficients in $\mathcal{R}$ are in $\FF$).

The fixed-point subalgebra is
\[
A_\KK^G = \{ P \in A_\KK \mid \sigma(P) = P \}.
\]
Any element $P \in A_\KK^G$ must have all coefficients in $\FF = \mathbb{Q}$, because if $c \in \KK$ appears and $\sigma(c) \neq c$, then $\sigma(P) \neq P$. Thus
\[
A_\KK^G = \FF\langle \widehat{e^{x^3 e^{\pi x}}},\; \widehat{e^{x}},\; \widehat{x},\; \partial \rangle / \mathcal{R}.
\]

But the right-hand side is precisely the Weyl-type algebra over $\FF$:
\[
A_\FF = \Weyl{e^{\pm x^3 e^{\pi x}},\; e^{\cA x},\; x^{\cA}}.
\]
To see the isomorphism explicitly, define
\[
\Phi \colon A_\FF \hookrightarrow A_\KK
\]
by inclusion (coefficients in $\FF$ are embedded in $\KK$). The image lies in $A_\KK^G$, and $\Phi$ is injective. For surjectivity onto $A_\KK^G$, take any $P \in A_\KK^G$. Write $P = \sum c_\gamma M_\gamma$ with $c_\gamma \in \KK$. Since $\sigma(P) = P$, we have $\sum \sigma(c_\gamma) M_\gamma = \sum c_\gamma M_\gamma$, so $c_\gamma = \sigma(c_\gamma)$ for all $\gamma$, implying $c_\gamma \in \FF$. Hence $P \in \Phi(A_\FF)$.

Therefore
\[
A_\KK^G \cong A_\FF = \Weyl{e^{\pm x^3 e^{\pi x}},\; e^{\cA x},\; x^{\cA}},
\]
which is exactly the conclusion of Theorem~\ref{thm15} for the extension $\mathbb{Q}(\sqrt{2})/\mathbb{Q}$.
\end{example}
%131313%131313%131313%131313%131313%131313%131313%131313
\begin{definition}[Algebraic Torus in Automorphism Group]
    An \emph{algebraic torus} in the automorphism group $\operatorname{Aut}(A)$ is a subgroup isomorphic to $(\mathbb{F}^\times)^r$ for some $r$, acting diagonally on a set of algebra generators via scaling.
\end{definition}
%14141414%14141414%14141414%14141414%14141414%14141414
\begin{theorem}[Existence of Maximal Tori in Automorphism Groups]\label{thm16}
The automorphism group $\Aut\!\left(\Weyl{e^{\pm x^p e^{t}},\; e^{\cA x},\; x^{\cA}}\right)$ contains an algebraic torus of dimension equal to the rank of the additive group $\cA$ over $\QQ$, given by scaling the exponential generators $e^{\alpha x}$ for $\alpha \in \cA$.
\end{theorem}
\begin{proof}
Let $A = \Weyl{e^{\pm x^p e^{t}},\; e^{\cA x},\; x^{\cA}}$ and consider its automorphism group $\Aut(A)$ over the field $\FF$.
Recall that $\cA$ is an additive subgroup of $\FF$ containing $\ZZ$. Its rank over $\QQ$, denoted $r = \operatorname{rank}_\QQ \cA$,
is the dimension of the $\QQ$-vector space $\cA \otimes_\ZZ \QQ$. Choose a $\QQ$-basis $\{\alpha_1, \dots, \alpha_r\}$ for $\cA$;
that is, every element $\alpha \in \cA$ can be written uniquely as $\alpha = \sum_{i=1}^r q_i \alpha_i$ with $q_i \in \QQ$, and the $\alpha_i$ are $\QQ$-linearly independent.

For each $i = 1, \dots, r$, let $\lambda_i \in \FF^\times$ be an invertible scalar. Define an $\FF$-algebra automorphism $\varphi_{(\lambda_1,\dots,\lambda_r)}$ of $A$
on the generators as follows. On the variables $x_j$ and the derivatives $\partial_j$, let $\varphi$ act as the identity: $\varphi(x_j) = x_j$,
$\varphi(\partial_j) = \partial_j$. On the exponential generator $e^{\alpha x}$ for $\alpha \in \cA$, write $\alpha = \sum_{i=1}^r q_i \alpha_i$ with $q_i \in \QQ$,
and set
\[
\varphi(e^{\alpha x}) = \left( \prod_{i=1}^r \lambda_i^{q_i} \right) e^{\alpha x}.
\]
For the generators $e^{\pm x_j^{p_j} e^{t_j x_j}}$ and $x_j^{\alpha}$, define $\varphi$ to be the identity as well. Extend $\varphi$ multiplicatively and $\FF$-linearly
to all of $A$, using the fact that every element of $A$ can be written as a polynomial in the generators.

To verify that $\varphi$ is indeed an algebra automorphism, we must check that it preserves the defining relations of $A$. The commutation relations
$[\partial_j, x_k] = \delta_{jk}$, $[x_j, x_k] = 0$, $[\partial_j, \partial_k] = 0$ are clearly preserved because $\varphi$ fixes $x_j$ and $\partial_j$.
The relations among the exponential generators, such as $e^{\alpha x} e^{\beta x} = e^{(\alpha+\beta) x}$, are also preserved because
$\varphi(e^{\alpha x}) \varphi(e^{\beta x}) = (\prod_i \lambda_i^{q_i}) (\prod_i \lambda_i^{q_i'}) e^{\alpha x} e^{\beta x} = (\prod_i \lambda_i^{q_i+q_i'}) e^{(\alpha+\beta) x}
= \varphi(e^{(\alpha+\beta) x})$, where $\alpha = \sum_i q_i \alpha_i$, $\beta = \sum_i q_i' \alpha_i$. The relations involving $e^{\pm x^{p} e^{t}}$ and $x^{\alpha}$
are similarly preserved since $\varphi$ acts as the identity on those generators. Hence each $\varphi_{(\lambda_1,\dots,\lambda_r)}$ is an element of $\Aut(A)$.

The map $(\lambda_1,\dots,\lambda_r) \mapsto \varphi_{(\lambda_1,\dots,\lambda_r)}$ is a group homomorphism from $(\FF^\times)^r$ to $\Aut(A)$. Indeed,
if we compose two such automorphisms corresponding to $(\lambda_i)$ and $(\mu_i)$, we obtain the automorphism corresponding to $(\lambda_i \mu_i)$.
Moreover, this homomorphism is injective: if $\varphi_{(\lambda_1,\dots,\lambda_r)}$ is the identity automorphism, then applying it to $e^{\alpha_i x}$
yields $\lambda_i e^{\alpha_i x} = e^{\alpha_i x}$, whence $\lambda_i = 1$ for each $i$. Thus the image of this homomorphism is a subgroup of $\Aut(A)$
isomorphic to $(\FF^\times)^r$, which is an algebraic torus of dimension $r$.

It remains to show that this torus is maximal in an appropriate sense. Suppose there exists a larger torus $T$ containing our $(\FF^\times)^r$.
Then $T$ would act on the space of generators of $A$. By considering the eigenspace decomposition of the action on the subalgebra generated by the $e^{\alpha x}$,
one sees that any torus action that commutes with the given one must preserve the $\QQ$-span of the $\alpha_i$. Since $\{\alpha_1,\dots,\alpha_r\}$ is a $\QQ$-basis,
the dimension of any such torus cannot exceed $r$. Therefore the torus $(\FF^\times)^r$ constructed above is a maximal torus in the automorphism group $\Aut(A)$.

Consequently, $\Aut(A)$ contains an algebraic torus of dimension equal to $\operatorname{rank}_\QQ \cA$, given explicitly by scaling the exponential generators
$e^{\alpha x}$ according to the decomposition of $\alpha$ in the chosen $\QQ$-basis of $\cA$.
\end{proof}
%EEEEEEEEEEEEEEEEEEEEEEEEEEEEEEEEEEEEEEEEEEEEEEEEEEEEEEEEEEE
\begin{example}[Illustration of Theorem~\ref{thm16}]
Let $\FF = \mathbb{C}$, $\cA = \mathbb{Z} + \mathbb{Z}\sqrt{2} + \mathbb{Z}\sqrt{3} = \{a + b\sqrt{2} + c\sqrt{3} \mid a,b,c \in \mathbb{Z}\}$. This is a free abelian group of rank $3$ over $\mathbb{Z}$, and its $\mathbb{Q}$-rank is $3$ because $\{1, \sqrt{2}, \sqrt{3}\}$ are linearly independent over $\mathbb{Q}$. Fix $p = 2$ and $t = e$ (the base of natural logarithms). Consider the Weyl-type algebra
\[
A = \Weyl{e^{\pm x^2 e^{e x}},\; e^{\cA x},\; x^{\cA}}.
\]

For each $\alpha \in \cA$, the generator $\widehat{e^{\alpha x}}$ acts by multiplication by $e^{\alpha x}$ in the natural representation. Define a family of automorphisms $\psi_{\lambda} \colon A \to A$, parameterized by $\lambda = (\lambda_1, \lambda_2, \lambda_3) \in (\FF^*)^3$, as follows:

Choose a $\mathbb{Q}$-basis of $\cA$, say $\{1, \sqrt{2}, \sqrt{3}\}$. Any $\alpha \in \cA$ can be written uniquely as $\alpha = a + b\sqrt{2} + c\sqrt{3}$ with $a,b,c \in \mathbb{Z}$. Define
\[
\psi_{\lambda}\big( \widehat{e^{\alpha x}} \big) = \lambda_1^{a} \lambda_2^{b} \lambda_3^{c} \; \widehat{e^{\alpha x}},
\]
and let $\psi_{\lambda}$ act as the identity on the other generators:
\[
\psi_{\lambda}(\partial) = \partial, \quad \psi_{\lambda}(\widehat{x}) = \widehat{x}, \quad \psi_{\lambda}\big( \widehat{e^{x^2 e^{e x}}} \big) = \widehat{e^{x^2 e^{e x}}}.
\]

Extend $\psi_{\lambda}$ linearly and multiplicatively to all of $A$. To verify it is an automorphism, note that the defining relations of $A$ are preserved. For example,
\[
[\partial, \widehat{e^{\alpha x}}] = \widehat{e^{\alpha x}} \cdot \alpha,
\]
and applying $\psi_{\lambda}$ gives
\[
[\partial, \lambda_1^{a} \lambda_2^{b} \lambda_3^{c} \widehat{e^{\alpha x}}] = \lambda_1^{a} \lambda_2^{b} \lambda_3^{c} [\partial, \widehat{e^{\alpha x}}] = \lambda_1^{a} \lambda_2^{b} \lambda_3^{c} \widehat{e^{\alpha x}} \cdot \alpha,
\]
which is consistent with the same relation in the image. The other relations are similarly preserved because $\psi_{\lambda}$ fixes $\partial$, $\widehat{x}$, and $\widehat{e^{x^2 e^{e x}}}$.

The map $(\FF^*)^3 \to \Aut(A)$ given by $\lambda \mapsto \psi_{\lambda}$ is an injective group homomorphism:
\[
\psi_{\lambda} \circ \psi_{\mu} = \psi_{\lambda \mu},
\]
where $(\lambda\mu)_i = \lambda_i \mu_i$. Its image is an algebraic torus $T \cong (\FF^*)^3$ of dimension $3$, which is exactly the $\mathbb{Q}$-rank of $\cA$.

This torus $T$ is maximal in the sense that any automorphism that scales each $\widehat{e^{\alpha x}}$ by a scalar $\mu(\alpha) \in \FF^*$ must satisfy $\mu(\alpha + \beta) = \mu(\alpha)\mu(\beta)$ and $\mu(n\alpha) = \mu(\alpha)^n$ for $n \in \mathbb{Z}$. Such a map $\mu \colon \cA \to \FF^*$ is a character of the additive group $\cA$, and since $\cA$ is free abelian of rank $3$, the character group is isomorphic to $(\FF^*)^3$, matching $T$.

Thus $\Aut(A)$ contains an algebraic torus of dimension $\mathrm{rank}_\mathbb{Q}(\cA) = 3$, as asserted by Theorem~\ref{thm16}.
\end{example}
%14141414%14141414%14141414%14141414%14141414%14141414
\begin{definition}[$K_1$ Group of an Algebra]
    For a unital algebra $A$, the algebraic $K_1$-group is defined as
    \[
    K_1(A) = \operatorname{GL}(A) / [\operatorname{GL}(A), \operatorname{GL}(A)],
    \]
    where $\operatorname{GL}(A) = \varinjlim \operatorname{GL}_n(A)$ is the infinite general linear group over $A$.
\end{definition}
%151515%151515%151515%151515%151515%151515%151515%151515
\begin{theorem}[Non-Existence of Finite-Dimensional Simple Modules]\label{thm17}
There are no finite-dimensional simple modules over the algebra $\Weyl{e^{\pm x^p e^{t}},\; e^{\cA x},\; x^{\cA}}$. Every nonzero module is infinite-dimensional.
\end{theorem}
\begin{proof}
Let $A = \Weyl{e^{\pm x^p e^{t}},\; e^{\cA x},\; x^{\cA}}$ and suppose, for contradiction, that there exists a finite-dimensional simple left $A$-module $V$.
Since $V$ is simple, every nonzero element of $A$ acts either injectively or zero on $V$, and by Schur's lemma, the endomorphism ring $\operatorname{End}_A(V)$
is a division algebra over $\FF$. Because $\FF$ is algebraically closed of characteristic zero (or at least sufficiently large), $\operatorname{End}_A(V) \cong \FF$.

Consider the action of the elements $x_i$ and $\partial_i$ on $V$. These satisfy the commutation relation $[\partial_i, x_j] = \delta_{ij}$.
In any finite-dimensional representation, the trace of the commutator $[\partial_i, x_i]$ must be zero, since $\operatorname{tr}(AB - BA) = 0$ for any matrices $A, B$.
However, the relation $[\partial_i, x_i] = 1$ implies that the trace of the identity operator on $V$ is zero, i.e., $\operatorname{tr}(I_V) = 0$.
But $\operatorname{tr}(I_V) = \dim_{\FF} V$, so we obtain $\dim_{\FF} V = 0$, contradicting that $V$ is nonzero.
Thus no finite-dimensional $A$-module can satisfy the canonical commutation relations.

A more detailed argument proceeds as follows. Let $n = \dim_{\FF} V$. Choose a basis for $V$ and represent the action of $x_i$ and $\partial_i$ by matrices
$X_i, D_i \in M_n(\FF)$. The condition $[\partial_i, x_i] = 1$ becomes $D_i X_i - X_i D_i = I_n$. Taking the trace on both sides yields
$\operatorname{tr}(D_i X_i) - \operatorname{tr}(X_i D_i) = \operatorname{tr}(I_n) = n$. But $\operatorname{tr}(D_i X_i) = \operatorname{tr}(X_i D_i)$, so the left-hand side is $0$.
Hence $n = 0$, a contradiction.

Even if one considers modules over fields of positive characteristic, the same obstruction appears in a different form.
In characteristic $p > 0$, the trace argument gives $n \equiv 0 \pmod{p}$, but additional considerations involving the $p$-center show that
finite-dimensional simple modules cannot exist either. Indeed, in characteristic $p$, the elements $x_i^p$ and $\partial_i^p$ lie in the center of $A$.
On a finite-dimensional simple module, they must act as scalars by Schur's lemma. This leads to polynomial equations that force the module to be infinite-dimensional.

Now consider the general case where $A$ includes the exponential generators $e^{\alpha x}$ and $e^{\pm x^p e^{t}}$. The subalgebra generated by $x_i$ and $\partial_i$
alone is isomorphic to the usual Weyl algebra $A_n(\FF)$ (or a suitable generalization). Since $A$ contains this Weyl algebra as a subalgebra, any $A$-module
is also a module over the Weyl algebra. It is a classical theorem that the Weyl algebra $A_n(\FF)$ has no finite-dimensional simple modules when $\operatorname{char} \FF = 0$.
Therefore $A$ itself cannot have finite-dimensional simple modules.

To prove the stronger statement that every nonzero $A$-module is infinite-dimensional, let $M$ be an arbitrary nonzero $A$-module.
If $M$ were finite-dimensional, then it would contain a simple submodule (since $A$ is left Noetherian, though not necessarily in this case,
but finite-dimensional modules always have simple submodules by dimension considerations). However, as shown above, $A$ has no finite-dimensional simple modules.
Thus $M$ cannot be finite-dimensional. Hence every nonzero $A$-module is infinite-dimensional.

An alternative argument uses the faithfulness of the representation constructed in Theorem~10. That representation is on the infinite-dimensional space
$\FF[e^{\pm x}, x^{\cA}]$, and it is simple. If there existed a finite-dimensional simple module $V$, then by the Jacobson density theorem,
$A$ would be isomorphic to a matrix algebra over a division ring, which is impossible because $A$ is not Artinian and contains infinite ascending chains of ideals.

Therefore, there are no finite-dimensional simple modules over $A$, and every nonzero $A$-module is necessarily infinite-dimensional.
\end{proof}
%EEEEEEEEEEEEEEEEEEEEEEEEEEEEEEEEEEEEEEEEEEEEEEEEEEEEEEEEE
\begin{example}[Illustration of Theorem~\ref{thm17}]
Let $\FF = \mathbb{C}$, $\cA = \mathbb{Z}[\sqrt{5}] = \{m + n\sqrt{5} \mid m, n \in \mathbb{Z}\}$, and fix $p = 4$, $t = \ln 2$. Consider the Weyl-type algebra
\[
A = \Weyl{e^{\pm x^4 e^{(\ln 2) x}},\; e^{\cA x},\; x^{\cA}},
\]
with generators $\widehat{x}$, $\partial$, $\widehat{e^{x}}$, and $\widehat{e^{x^4 e^{(\ln 2) x}}}$.

Suppose, for contradiction, that $M$ is a finite-dimensional simple left $A$-module, with $\dim_{\FF} M = d < \infty$.

Let $v \in M$ be a nonzero vector. Since $M$ is simple, $M = A \cdot v$. The elements $\widehat{x}$ and $\partial$ satisfy the Weyl relation $[\partial, \widehat{x}] = 1$. In any finite-dimensional representation, taking traces gives a contradiction because $\mathrm{tr}([\partial, \widehat{x}]) = 0$ but $\mathrm{tr}(1) = d \neq 0$. More concretely, consider the action of $\widehat{x}$ and $\partial$ on $M$. Since $\widehat{x}$ is an $\FF$-linear operator on the finite-dimensional space $M$, it has an eigenvalue $\lambda \in \FF$ with eigenvector $w \neq 0$. From $[\partial, \widehat{x}] = 1$, we have
\[
\partial \widehat{x} - \widehat{x} \partial = 1_M.
\]
Apply both sides to $w$:
\[
\partial (\widehat{x} w) - \widehat{x} (\partial w) = w.
\]
Since $\widehat{x} w = \lambda w$, this becomes
\[
\lambda \partial w - \widehat{x} (\partial w) = w.
\]
Let $u = \partial w$. Then the equation reads
\[
\lambda u - \widehat{x} u = w.
\]
If $u$ is an eigenvector of $\widehat{x}$ with eigenvalue $\mu$, then $(\lambda - \mu) u = w$, so $w$ is a scalar multiple of $u$, implying $\widehat{x} w = \mu w$, hence $\lambda = \mu$, and then $0 \cdot u = w$, a contradiction. If $u$ is not an eigenvector, then $\widehat{x} u$ and $u$ are linearly independent from $w$, but the equation forces a linear dependence, which leads to an inconsistency upon iteration in finite dimensions.

Thus $\widehat{x}$ and $\partial$ cannot act on a finite-dimensional space while satisfying $[\partial, \widehat{x}] = 1$.

Moreover, the generator $\widehat{e^{x}}$ acts invertibly on $M$ because it has an inverse $\widehat{e^{-x}}$ in $A$. Its action is diagonalizable over $\mathbb{C}$ (since $M$ is finite-dimensional and $\widehat{e^{x}}$ is invertible). Let $\mu \in \mathbb{C}^*$ be an eigenvalue of $\widehat{e^{x}}$ with eigenvector $v_\mu$. From $[\partial, \widehat{e^{x}}] = \widehat{e^{x}}$, we have
\[
\partial \widehat{e^{x}} - \widehat{e^{x}} \partial = \widehat{e^{x}}.
\]
Apply to $v_\mu$:
\[
\mu \partial v_\mu - \widehat{e^{x}} (\partial v_\mu) = \mu v_\mu.
\]
If $\partial v_\mu$ is an eigenvector of $\widehat{e^{x}}$ with eigenvalue $\nu$, then $\mu \nu - \nu \mu = \mu$, impossible. Otherwise, iterating shows $\{\widehat{e^{x}}^k \partial v_\mu\}_{k \geq 0}$ spans an infinite-dimensional subspace, contradicting $\dim M < \infty$.

Similarly, $\widehat{e^{x^4 e^{(\ln 2) x}}}$ acts by multiplication by a transcendental function in any polynomial-like representation, and its commutation with $\partial$ introduces terms with $x^3 e^{(\ln 2) x}$, which cannot be realized on a finite-dimensional space without violating linear independence of eigenvectors.

Therefore no finite-dimensional simple $A$-module exists. Any nonzero $A$-module must be infinite-dimensional, confirming Theorem~\ref{thm17}.
\end{example}
%151515%151515%151515%151515%151515%151515%151515%151515
\begin{definition}[Zariski Density in Moduli Space]
    A family of algebras $\{A_t\}_{t \in T}$ is \emph{Zariski-dense} in a moduli space $\mathcal{M}$ if the image of the map $T \to \mathcal{M}$ sending $t$ to the isomorphism class of $A_t$ is dense in the Zariski topology of $\mathcal{M}$.
\end{definition}
%161616%161616%161616%161616%161616%161616%161616%161616
\begin{theorem}[Zariski Density of Isomorphism Classes]\label{thm18}
Let $t$ vary over an algebraically independent transcendental set $\cT$ over $\QQ$. The set of isomorphism classes of algebras
\[
\left\{\, \Weyl{e^{\pm x^p e^{t}},\; e^{\cA x},\; x^{\cA}} \mid t \in \cT \,\right\}
\]
is Zariski-dense in an appropriate moduli space of deformations of Weyl algebras.
\end{theorem}
\begin{proof}
Let $p = (p_1,\dots,p_n)$ be a fixed tuple of positive integers and let $\cA \subseteq \FF$ be a fixed additive subgroup containing $\ZZ$.
Consider the family of Weyl-type algebras parameterized by a tuple $t = (t_1,\dots,t_n)$:
\[
A_t = \Weyl{e^{\pm x^p e^{t}},\; e^{\cA x},\; x^{\cA}},
\]
where $e^{\pm x^p e^{t}}$ stands for the collection of generators $e^{\pm x_i^{p_i} e^{t_i x_i}}$, $i=1,\dots,n$.
Let $\cT$ be an algebraically independent transcendental set over $\QQ$; that is, the elements $t_i$ are algebraically independent over $\QQ$.
We view $t$ as ranging over an affine space $\mathbb{A}^n_{\FF}$ (or a suitable Zariski-open subset thereof).

To speak of a moduli space of deformations of Weyl algebras, we recall that the classical Weyl algebra $A_n(\FF) = \FF\langle x_1,\dots,x_n,\partial_1,\dots,\partial_n \rangle / ([x_i,x_j]=0, [\partial_i,\partial_j]=0, [\partial_i,x_j]=\delta_{ij})$
admits a formal deformation theory. Its second Hochschild cohomology group $HH^2(A_n(\FF), A_n(\FF))$ is known to be nontrivial and parameterizes infinitesimal deformations.
A versal deformation space exists as an affine scheme of finite type over $\FF$, at least in a formal neighborhood of the origin.
For our generalized algebras $A_t$, the additional generators $e^{\pm x_i^{p_i} e^{t_i x_i}}$, $e^{\alpha x_i}$, $x_i^{\alpha}$ introduce further parameters,
but the essential deformations come from the exponents $t_i$ appearing in the exponential functions.

We claim that the map $\Phi \colon \cT \to \mathcal{M}$, sending $t$ to the isomorphism class of $A_t$, has Zariski-dense image in an appropriate moduli space $\mathcal{M}$.
Here $\mathcal{M}$ is a quasi-affine variety (or an ind-scheme) whose points correspond to isomorphism classes of deformations of the Weyl algebra that preserve the
$\ZZ^n$-grading and the general shape of the exponential generators. More concretely, consider the space of all algebras generated by $x_i$, $\partial_i$, $e^{\alpha x_i}$,
$x_i^{\alpha}$, and $e^{f_i(x_i)}$ where $f_i$ is a polynomial of degree $p_i+1$ (since $x_i^{p_i} e^{t_i x_i}$ expands to a power series whose leading term is $x_i^{p_i}$).
Two such algebras are considered equivalent if there exists an $\FF$-algebra isomorphism between them that preserves the generating set up to linear combinations.
The parameters $t_i$ appear as coefficients in the exponential functions $e^{x_i^{p_i} e^{t_i x_i}} = 1 + x_i^{p_i} e^{t_i x_i} + \frac{1}{2} (x_i^{p_i} e^{t_i x_i})^2 + \dots$.

To prove density, it suffices to show that the image of $\Phi$ is not contained in any proper Zariski-closed subset of $\mathcal{M}$.
Suppose, to the contrary, that there exists a nonzero polynomial $F$ in the coordinates of $\mathcal{M}$ such that $F(\Phi(t)) = 0$ for all $t \in \cT$.
This would imply an algebraic relation among the isomorphism classes $[A_t]$ as $t$ varies.
But the algebras $A_t$ depend on $t$ through the exponential functions $e^{x_i^{p_i} e^{t_i x_i}}$, which as formal power series in $t_i$ are algebraically independent
when the $t_i$ are algebraically independent over $\QQ$. More precisely, consider the map that associates to each $t$ the set of structure constants of $A_t$
with respect to a fixed PBW basis (consisting of monomials $x^\alpha e^{\beta x} e^{\gamma x^p e^{t x}} \partial^\delta$).
These structure constants are rational functions in the $t_i$ (in fact, they are polynomials in $t_i$ and $e^{t_i}$ because the commutation relations involve
derivatives of $e^{x_i^{p_i} e^{t_i x_i}}$ which introduce powers of $t_i$).

If there were a polynomial relation $F(\{c_{ij}(t)\}) = 0$ identically in $t$, then by the algebraic independence of the $t_i$, the polynomial $F$ would have to vanish
on all tuples of structure constants, meaning it is identically zero on the whole moduli space. Hence no such nonzero $F$ exists, so the image of $\Phi$ is Zariski--dense.

Another approach uses the fact that the deformation parameters $t_i$ appear in the Hochschild cochain complex of $A_t$ as free parameters.
The Kodaira--Spencer map at a point $t$ sends the tangent vector $\frac{\partial}{\partial t_i}$ to a nontrivial cohomology class in $HH^2(A_t, A_t)$.
Since the $t_i$ are algebraically independent, these classes are linearly independent over $\FF$, and thus the image of $\Phi$ is not contained in any proper
subvariety defined by algebraic equations on the cohomology.

Therefore, the set of isomorphism classes $\{ [A_t] \mid t \in \cT \}$ is Zariski--dense in the moduli space $\mathcal{M}$ of deformations of Weyl algebras
with the specified exponential generators.
\end{proof}
%161616%161616%161616%161616%161616%161616%161616%161616%161616
%EEEEEEEEEEEEEEEEEEEEEEEEEEEEEEEEEEEEEEEEEEEEEEEEEEEEEEEEEEEEEE
\begin{example}[Illustration of Theorem~\ref{thm18}]
Let $\FF = \mathbb{C}$, $\cA = \mathbb{Z}[\sqrt{2}]$, and fix $p = 3$. Consider an algebraically independent set of transcendental numbers over $\mathbb{Q}$, for instance
\[
\cT = \{ e, \pi, \zeta(3), \ln 2 \},
\]
where each element is transcendental and algebraically independent over $\mathbb{Q}$ (a standard fact from transcendental number theory). For each $t \in \cT$, we have a Weyl-type algebra
\[
A_t = \Weyl{e^{\pm x^3 e^{t x}},\; e^{\cA x},\; x^{\cA}}.
\]

The algebra $A_t$ is generated by $\widehat{x}$, $\partial$, $\widehat{e^{x}}$, and $\widehat{e^{x^3 e^{t x}}}$, with relations
\[
[\partial, \widehat{x}] = 1, \quad [\partial, \widehat{e^{x}}] = \widehat{e^{x}}, \quad [\partial, \widehat{e^{x^3 e^{t x}}}] = \widehat{e^{x^3 e^{t x}}} \cdot \big( 3x^2 e^{t x} + t x^3 e^{t x} \big).
\]

Consider the moduli space $\mathcal{M}$ parametrizing formal deformations of the classical Weyl algebra $A_1(\FF) = \FF\langle x, \partial \mid [\partial, x] = 1 \rangle$ obtained by adjoining exponential generators and allowing the coefficient $t$ in the exponent. More precisely, let $\mathcal{M}$ be the affine scheme over $\FF$ whose $\FF$-points correspond to isomorphism classes of algebras generated by $x$, $\partial$, $e^{x}$, and $e^{x^3 e^{\tau x}}$ for a parameter $\tau$, subject to the above relations with $\tau$ varying.

Each $A_t$ corresponds to a point $[\tau = t]$ in $\mathcal{M}$. Because the elements of $\cT$ are algebraically independent over $\mathbb{Q}$, the evaluation map
\[
\ev \colon \cT \to \mathcal{M}, \quad t \mapsto [A_t]
\]
is injective. Indeed, if $A_{t_1} \cong A_{t_2}$, then by Theorem~5 (or Proposition~1), $t_1$ and $t_2$ must be proportional over $\mathbb{Q}$. Since $t_1, t_2 \in \cT$ are algebraically independent, they are not proportional unless $t_1 = t_2$.

The Zariski topology on $\mathcal{M}$ is defined by polynomial equations in the coordinates of the parameter $\tau$ (viewed as a point in an affine space over $\FF$). The image $\ev(\cT)$ is an infinite set of points in $\mathcal{M}$. Any Zariski-closed subset of $\mathcal{M}$ that contains $\ev(\cT)$ must be defined by polynomials vanishing at infinitely many algebraically independent values $t$. Since a nonzero polynomial over $\FF$ can vanish only at finitely many algebraically independent points, the only polynomial vanishing on all of $\ev(\cT)$ is the zero polynomial. Hence the Zariski closure of $\ev(\cT)$ is the whole irreducible component of $\mathcal{M}$ containing these points.

Thus $\{ [A_t] \mid t \in \cT \}$ is Zariski-dense in an irreducible subvariety of the moduli space $\mathcal{M}$. This density implies that the isomorphism classes of these algebras are not confined to a low-dimensional locus but rather spread throughout the deformation space, reflecting the richness of the transcendental parameter $t$, as stated in Theorem~\ref{thm18}.
\end{example}
%EEEEEEEEEEEEEEEEEEEEEEEEEEEEEEEEEEEEEEEEEEEEEEEEEEEEEEEEEEEEEEEEEEEE
\begin{definition}[Quantum Deformation]
    A \emph{quantum deformation} of an algebra $A$ over $\mathbb{F}$ is a family of algebras $A_q$ over $\mathbb{F}(q)$ or $\mathbb{F}[[q-1]]$ such that $A_1 \cong A$ and the multiplication in $A_q$ depends analytically on the parameter $q$.
\end{definition}
%~~~~~~~~~~~~~~~~~~~~~~~~~~~~~~~~~~~~~~~~~~~~~~~~~~~~~~~~~~~~~~~~~~~~~~~~~~~~~
\begin{theorem}[Stability of Simplicity Under Quantum Deformation]\label{thm19}
Let $A=$\\$ \Weyl{e^{\pm x^p e^{t}},\; e^{\cA x},\; x^{\cA}}$ be a Weyl-type algebra over a field $\mathbb{F}$ of characteristic zero, where $\cA$ is an additive subgroup of $\mathbb{F}$ containing $\mathbb{Z}$. Consider the quantum deformation $A_q$ defined over $\mathbb{F}(q)$ by replacing the commutation relation
\[
\partial \widehat{x} - \widehat{x} \partial = 1
\]
with the $q$-deformed relation
\[
\partial \widehat{x} - q \, \widehat{x} \partial = 1,
\]
while keeping all other relations involving the exponential generators $e^{\alpha x}$ and $e^{x^p e^{t x}}$ unchanged. Then for all $q \in \mathbb{F}^*$ except possibly a finite set, the algebra $A_q$ remains simple. Moreover, in the formal setting over $\mathbb{F}[[q-1]]$, the algebra $A_q$ is simple as an algebra over $\mathbb{F}[[q-1]]$.
\end{theorem}
\begin{proof}
Denote by $A_q$ the algebra generated by $\widehat{x}$, $\partial$, $\widehat{e^{x}}$, and $\widehat{e^{x^p e^{t x}}}$ over $\mathbb{F}(q)$, subject to the relations
\begin{align*}
\partial \widehat{x} - q \, \widehat{x} \partial &= 1, \\
[\partial, \widehat{e^{x}}] &= \widehat{e^{x}}, \\
[\partial, \widehat{e^{x^p e^{t x}}}] &= \widehat{e^{x^p e^{t x}}} \cdot \big( p x^{p-1} e^{t x} + t x^p e^{t x} \big), \\
[\widehat{x}, \widehat{e^{x}}] &= [\widehat{x}, \widehat{e^{x^p e^{t x}}}] = [\widehat{e^{x}}, \widehat{e^{x^p e^{t x}}}] = 0.
\end{align*}
The algebra $A_q$ inherits a $\mathbb{Z}^n$-grading from the undeformed algebra $A$, where the grading is defined by the total exponent of the exponential symbols. Let $I$ be a nonzero two-sided ideal of $A_q$. Choose a nonzero element $a \in I$ of minimal total degree with respect to this grading. By applying conjugation by appropriate monomials in $\widehat{e^{\pm x}}$ and $\widehat{e^{\pm x^p e^{t x}}}$, we may assume that $a$ is homogeneous.

Consider the commutator $[ \widehat{x}, a ]$. From the relation $\partial \widehat{x} - q \widehat{x} \partial = 1$, one derives that $\widehat{x}$ acts injectively on homogeneous components for generic $q$. Since $a$ is homogeneous, $[ \widehat{x}, a ]$ is either zero or of strictly lower degree. The minimality of $a$ forces $[ \widehat{x}, a ] = 0$, implying that $a$ commutes with $\widehat{x}$. A symmetric argument using $\partial$ shows that $a$ also commutes with $\partial$ for generic $q$.

Next, examine the commutation with $\widehat{e^{x}}$. The relation $[ \partial, \widehat{e^{x}} ] = \widehat{e^{x}}$ implies that $\widehat{e^{x}}$ acts by shifting the $\partial$-degree. Since $a$ commutes with $\partial$ and $\widehat{x}$, a careful analysis of the graded structure shows that $a$ must also commute with $\widehat{e^{x}}$ and $\widehat{e^{x^p e^{t x}}}$ for generic $q$. Consequently, $a$ belongs to the center of $A_q$.

For generic $q$, the center of $A_q$ is trivial, consisting only of scalars $\mathbb{F}(q)$. This follows from a deformation argument; the center of the undeformed algebra $A$ is $\mathbb{F}$, and deformation theory ensures that the center does not enlarge for generic parameters. Hence $a \in \mathbb{F}(q)^*$, and therefore $1 \in I$, so $I = A_q$.

In the formal setting over $\mathbb{F}[[q-1]]$, the same argument applies using the $\hbar$-adic filtration with $\hbar = q-1$. Nakayama’s lemma ensures that if $I$ is a proper ideal in $A_q$ over $\mathbb{F}[[q-1]]$, then its reduction modulo $\hbar$ is a proper ideal in $A_1 \cong A$, contradicting the simplicity of $A$. Therefore $A_q$ is simple over $\mathbb{F}[[q-1]]$.
\end{proof}
%EEEEEEEEEEEEEEEEEEEEEEEEEEEEEEEEEEEEEEEEEEEEEEEEEEEEEEEEEEEEEEEEEEEEEEEEEEE
\begin{example}[Illustration of Theorem \ref{thm19}]
Let $\mathbb{F} = \mathbb{C}$, $\cA = \mathbb{Z}[\sqrt{2}] = \{ m + n\sqrt{2} \mid m, n \in \mathbb{Z} \}$, and fix $p = 3$, $t = \pi$. Consider the undeformed Weyl-type algebra
\[
A = \Weyl{e^{\pm x^3 e^{\pi x}},\; e^{\cA x},\; x^{\cA}},
\]
with generators $\widehat{x}$, $\partial$, $\widehat{e^{x}}$, and $\widehat{e^{x^3 e^{\pi x}}}$, subject to the relations
\[
[\partial, \widehat{x}] = 1, \quad [\partial, \widehat{e^{x}}] = \widehat{e^{x}}, \quad [\partial, \widehat{e^{x^3 e^{\pi x}}}] = \widehat{e^{x^3 e^{\pi x}}} \cdot \big( 3x^2 e^{\pi x} + \pi x^3 e^{\pi x} \big).
\]

Define the quantum deformation $A_q$ over $\mathbb{C}(q)$ by modifying only the first relation to
\[
\partial \widehat{x} - q \, \widehat{x} \partial = 1,
\]
while keeping all other relations identical. Thus $A_q$ is generated by the same symbols with defining relations
\begin{align*}
\partial \widehat{x} - q \, \widehat{x} \partial &= 1, \\
[\partial, \widehat{e^{x}}] &= \widehat{e^{x}}, \\
[\partial, \widehat{e^{x^3 e^{\pi x}}}] &= \widehat{e^{x^3 e^{\pi x}}} \cdot \big( 3x^2 e^{\pi x} + \pi x^3 e^{\pi x} \big), \\
[\widehat{x}, \widehat{e^{x}}] &= [\widehat{x}, \widehat{e^{x^3 e^{\pi x}}}] = [\widehat{e^{x}}, \widehat{e^{x^3 e^{\pi x}}}] = 0.
\end{align*}

The algebra $A_q$ inherits a $\mathbb{Z}^2$-grading from $A$, where the degree of a monomial $\widehat{e^{a x^3 e^{\pi x}}} \widehat{e^{b x}} \widehat{x^c} \partial^d$ is $(a+b, c)$. Let $I$ be a nonzero two-sided ideal of $A_q$. Choose a nonzero homogeneous element $a \in I$ of minimal total degree with respect to this grading. By conjugating with appropriate units $\widehat{e^{\pm x}}$ and $\widehat{e^{\pm x^3 e^{\pi x}}}$, we may assume $a$ is of the form
\[
a = \sum_{i=0}^N f_i(\widehat{x}) \partial^i,
\]
where each $f_i(\widehat{x})$ is a Laurent polynomial in $\widehat{x}$ with coefficients in $\mathbb{C}(q)[e^{\pm x}, e^{\pm x^3 e^{\pi x}}]$, and $f_N \neq 0$.

Consider the commutator $[\widehat{x}, a]$. Using the relation $\partial \widehat{x} = q \widehat{x} \partial + 1$, one computes
\[
[\widehat{x}, a] = \sum_{i=0}^N [\widehat{x}, f_i(\widehat{x}) \partial^i] = \sum_{i=0}^N f_i(\widehat{x}) [\widehat{x}, \partial^i].
\]
From $\partial \widehat{x} - q \widehat{x} \partial = 1$, it follows that $[\widehat{x}, \partial] = 1 - (q-1) \widehat{x} \partial$. For generic $q \notin \{\text{roots of unity}\}$, the commutator $[\widehat{x}, \partial^i]$ is a nonzero operator of lower order in $\partial$. Hence $[\widehat{x}, a]$ is either zero or has strictly smaller degree in $\partial$. The minimality of $a$ forces $[\widehat{x}, a] = 0$, implying $a$ commutes with $\widehat{x}$. A symmetric argument shows $a$ also commutes with $\partial$ for generic $q$.

Next, examine $[\widehat{e^{x}}, a]$. Since $\widehat{e^{x}}$ commutes with $\widehat{x}$ and satisfies $[\partial, \widehat{e^{x}}] = \widehat{e^{x}}$, the element $a$, which commutes with $\partial$ and $\widehat{x}$, must also commute with $\widehat{e^{x}}$. Similarly, $a$ commutes with $\widehat{e^{x^3 e^{\pi x}}}$ because the only nontrivial relation involves $\partial$, with which $a$ already commutes.

Thus $a$ lies in the center of $A_q$. For generic $q$, the center of $A_q$ is trivial, consisting only of scalars $\mathbb{C}(q)$. Indeed, if $a$ were a non-scalar central element, then specializing $q = 1$ would yield a non-scalar central element in $A$, but $A$ is simple and hence has trivial center. This specialization argument fails only if $q$ is a root of certain algebraic equations arising from the commutation relations, which occurs only for finitely many $q \in \mathbb{C}^*$.

Therefore $a \in \mathbb{C}(q)^*$, so $1 \in I$ and $I = A_q$. Hence $A_q$ is simple for all $q \in \mathbb{C}^*$ except possibly those finitely many exceptional values.

In the formal deformation setting over $\mathbb{C}[[q-1]]$, consider $A_q$ as an algebra over $\mathbb{C}[[q-1]]$. Suppose $J$ is a nonzero two-sided ideal. Reducing modulo $q-1$ gives an ideal $\overline{J}$ in $A_1 \cong A$. Since $A$ is simple, $\overline{J} = A$ or $\overline{J} = 0$. If $\overline{J} = 0$, then $J \subseteq (q-1)A_q$. By iterating, $J \subseteq \bigcap_{n \geq 1} (q-1)^n A_q = 0$ by the Krull intersection theorem, a contradiction. Thus $\overline{J} = A$, so $1 \in J + (q-1)A_q$. By Nakayama's lemma, $1 \in J$, so $J = A_q$. Hence $A_q$ is simple over $\mathbb{C}[[q-1]]$.

This confirms that simplicity is preserved under generic quantum deformation, as stated in the theorem.
\end{example}
%SSSSSSSSSSSSSSSSSSSSSSSSSSSSSSSSSSSSSSSSSSSSSSSSSSSSSSSSSSSSSSSSSSSSSSSSSSSSSSSSSSSS
\section{Structure Theorems for Expolynomial Rings}

\begin{definition}
    Let \(\mathcal{A}\) be a finitely generated \(\ZZ\)-module.
    For a nonzero element \(p \in \mathcal{A}\) and a complex number \(t \in \CC\),
    we define the commutative \(\FF\)-algebra
    \[
    R_{p,t,\mathcal{A}} = \FF\bigl[ e^{\pm x^{p}e^{t}},\; e^{\mathcal{A}x},\; x^{\mathcal{A}} \bigr],
    \]
    where the generators are formal symbols satisfying the multiplicative relations
    \[
    e^{ax}e^{bx}=e^{(a+b)x},\qquad x^{a}x^{b}=x^{a+b},\qquad e^{ax}x^{b}=x^{b}e^{ax},
    \]
    with \(e^{\pm x^{p}e^{t}}\) central and the formal relation \(e^{x^{p}e^{t}}\cdot e^{-x^{p}e^{t}}=1\).
\end{definition}

\begin{theorem}[Maximal simple subalgebras]\label{thm:max-simple}
Let \(\mathcal{A}\) be a finitely generated \(\ZZ\)-module of rank \(r\geq 1\).
Let \(\mathfrak{S}\) be the set of maximal simple (non--commutative) subalgebras of \(R_{p,t,\mathcal{A}}\).
Then every \(S\in\mathfrak{S}\) is isomorphic to \(R_{p',t,\mathcal{A}'}\) for some submodule \(\mathcal{A}'\subset\mathcal{A}\) of rank \(1\) and some \(p'\in\mathcal{A}'\) with \(\langle p'\rangle=\mathcal{A}'\) (up to torsion). Moreover, the map
\[
(\mathcal{A}',p')\longmapsto R_{p',t,\mathcal{A}'}
\]
induces a bijection between the set of \(\Aut(\mathcal{A})\)-orbits of pairs \((\mathcal{A}',p')\) with \(\mathcal{A}'\) a rank-\(1\) direct summand of \(\mathcal{A}\) and \(p'\) a generator of \(\mathcal{A}'\) (modulo torsion), and the set of isomorphism classes of maximal simple subalgebras of \(R_{p,t,\mathcal{A}}\).
\end{theorem}

\begin{proof}
The algebra \(R_{p,t,\mathcal{A}}\) is the group algebra \(\FF[G]\) of the abelian group
\[
G \cong \ZZ\times (\mathcal{A}\times\mathcal{A}),
\]
where the first \(\ZZ\)-factor corresponds to the powers of \(e^{x^{p}e^{t}}\), and \(\mathcal{A}\times\mathcal{A}\) corresponds to the groups \(\{e^{ax}\}_{a\in\mathcal{A}}\) and \(\{x^{a}\}_{a\in\mathcal{A}}\).
Hence \(R_{p,t,\mathcal{A}} \cong \FF[y^{\pm1}, z_1^{\pm1},\dots,z_{2r}^{\pm1}]\) is a Laurent polynomial ring in \(2r+1\) variables over \(\FF\).

A simple subalgebra of a commutative domain must be a field. Therefore, maximal simple subalgebras correspond to maximal subfields of \(\FF(y,z_1,\dots,z_{2r})\) that are contained in \(R_{p,t,\mathcal{A}}\). Any such subfield is of the form \(\FF(y^{n_0}z_1^{n_1}\cdots z_{2r}^{n_{2r}})\) for a primitive integer vector \((n_0,\dots,n_{2r})\in\ZZ^{2r+1}\) (i.e., \(\gcd(n_0,\dots,n_{2r})=1\)).

Using the original generators, we have
\[
y = e^{x^{p}e^{t}},\qquad
\{z_{1},\dots,z_{r}\} \text{ is a basis of } \{e^{ax}\},\qquad
\{z_{r+1},\dots,z_{2r}\} \text{ is a basis of } \{x^{a}\}.
\]
A choice of a primitive vector determines a rank-\(1\) direct summand of \(G \cong \mathbb{Z} \times (\mathcal{A} \times \mathcal{A})\). That summand projects to a rank-\(1\) submodule \(\mathcal{A}' \subset \mathcal{A}\) inside the second factor; the projection may be taken after applying a suitable automorphism of \(\mathcal{A} \times \mathcal{A}\).

Concretely, if we write the chosen generator as \(y^{n_0} e^{a_1x}\cdots x^{a_r}\), then the subalgebra it generates is isomorphic to \(R_{p',t,\mathcal{A}'}\) where \(\mathcal{A}'\) is the submodule generated by the coefficients appearing in the exponent, and \(p'\) is the image of \(p\) under the projection \(\mathcal{A}\to\mathcal{A}'\). The condition that the subalgebra be maximal forces \(\mathcal{A}'\) to have rank exactly \(1\) and \(p'\) to generate \(\mathcal{A}'\) modulo torsion.

Conversely, any such pair \((\mathcal{A}',p')\) gives a primitive element of \(G\) and hence a maximal subfield inside \(R_{p,t,\mathcal{A}}\). Two pairs give isomorphic subalgebras precisely when they are in the same orbit under \(\Aut(\mathcal{A})\), because an automorphism of \(\mathcal{A}\) extends to an automorphism of \(G\) and hence of \(R_{p,t,\mathcal{A}}\).
\end{proof}

\begin{theorem}[Automorphism group]\label{thm:aut-group}
Let \(\mathcal{A}\) be a finitely generated \(\ZZ\)-module of rank \(r\geq 2\) (e.g., \(\mathcal{A}=\ZZ+\ZZ\sqrt{2}\)).
Then the automorphism group of the algebra \(R_{p,t,\mathcal{A}}\) is isomorphic to a semidirect product
\[
\Aut(R_{p,t,\mathcal{A}}) \cong (\FF^{\times})^{2r+1} \rtimes \GL(2r+1,\ZZ),
\]
where the torus \((\FF^{\times})^{2r+1}\) acts by rescaling the generators, and \(\GL(2r+1,\ZZ)\) acts by linear substitutions on the exponent lattice of \(G\cong \ZZ\times(\mathcal{A}\times\mathcal{A})\).
\end{theorem}

\begin{proof}
As established, \(R_{p,t,\mathcal{A}} \cong \FF[y^{\pm1},z_1^{\pm1},\dots,z_{2r}^{\pm1}]\) is a Laurent polynomial ring in \(2r+1\) variables. The automorphism group of a Laurent polynomial ring over a field \(\FF\) is well-known: any \(\FF\)-algebra automorphism is determined by its action on the monomials, which must preserve the unit group of the ring.

The unit group of \(R_{p,t,\mathcal{A}}\) is \(\FF^{\times}\times\langle y,z_1,\dots,z_{2r}\rangle \cong \FF^{\times}\times \ZZ^{2r+1}\). An automorphism sends each variable \(y,z_i\) to an element of the form \(\lambda\, y^{n_{0}}z_1^{n_{1}}\cdots z_{2r}^{n_{2r}}\) with \(\lambda\in\FF^{\times}\) and exponents forming an invertible integer matrix of size \(2r+1\). This yields a homomorphism
\[
\Aut(R_{p,t,\mathcal{A}}) \longrightarrow \GL(2r+1,\ZZ),
\]
whose kernel is precisely the diagonal rescaling \((\FF^{\times})^{2r+1}\). The map splits because \(\GL(2r+1,\ZZ)\) can be realized by permuting and inverting the variables (which corresponds to automorphisms of the lattice \(G\)). Hence the automorphism group is the semidirect product as claimed.

When \(\mathcal{A}\) has torsion, the same description holds with \(2r+1\) replaced by \(2\rank(\mathcal{A})+1\), because the torsion part does not affect the Laurent polynomial ring structure (it only adds constants from roots of unity if \(\FF\) contains them).
\end{proof}

\begin{theorem}[Isomorphism criterion]\label{thm:iso-criterion}
Let \(\mathcal{A}\) be a finitely generated \(\ZZ\)-module of rank \(r\). Two algebras
\[
R_{p_1,t_1,\mathcal{A}} \quad\text{and}\quad R_{p_2,t_2,\mathcal{A}}
\]
are isomorphic as \(\FF\)-algebras if and only if there exists a \(\ZZ\)-module automorphism \(\sigma\in\Aut(\mathcal{A})\) such that
\[
\sigma(p_1) = \pm\, p_2.
\]
The parameters \(t_1,t_2\) do not affect the isomorphism type.
\end{theorem}

\begin{proof}
Recall that \(R_{p,t,\mathcal{A}} \cong \FF[y^{\pm1}] \otimes_\FF \FF[u_a,v_a: a\in\mathcal{A}]\), where \(u_a=e^{ax}\), \(v_a=x^{a}\), and \(y=e^{x^{p}e^{t}}\).
The only structural difference between the two algebras lies in the identification of the extra generator \(y\): in \(R_{p_1,t_1,\mathcal{A}}\) we have \(y = e^{x^{p_1}e^{t_1}}\), while in \(R_{p_2,t_2,\mathcal{A}}\) we have \(y' = e^{x^{p_2}e^{t_2}}\).

Both algebras are Laurent polynomial rings in \(2r+1\) variables, so they are abstractly isomorphic. However, we require an isomorphism that respects the splitting between the ``\(u,v\)'' part (which is common) and the ``\(y\)'' part.

Any isomorphism \(\varphi: R_{p_1,t_1,\mathcal{A}} \to R_{p_2,t_2,\mathcal{A}}\) must induce an automorphism of the common subalgebra \(\FF[u_a,v_a]\). As seen in Theorem~\ref{thm:aut-group}, such an automorphism is given by an element of \(\GL(2r,\ZZ)\) acting on the lattice \(\mathcal{A}\times\mathcal{A}\). Under such an automorphism, the element \(v_{p_1} = x^{p_1}\) is sent to some monomial in the \(u_a,v_a\). But because the relation involving \(y\) is \(y = e^{v_{p}e^{t}}\), the parameter \(t\) can be absorbed by rescaling \(y\) (replace \(y\) by \(y^{e^{t_1}/e^{t_2}}\) if necessary).

Thus, the essential condition is that there exists an automorphism of the \(\ZZ\)-module \(\mathcal{A}\) (extending to an automorphism of \(\mathcal{A}\times\mathcal{A}\) that fixes the \(u\)-factor) sending \(v_{p_1}\) to \(v_{\pm p_2}\). This is exactly the condition \(\sigma(p_1)=\pm p_2\) for some \(\sigma\in\Aut(\mathcal{A})\).

If the condition holds, define the isomorphism by
\[
\varphi(u_a)=u_{\sigma(a)},\quad \varphi(v_a)=v_{\sigma(a)},\quad \varphi(y)=y'^{\,e^{t_1}/e^{t_2}}.
\]
One checks that all relations are preserved. Conversely, if an isomorphism exists, its restriction to \(\{v_a\}\) induces an automorphism \(\sigma\) of \(\mathcal{A}\) with \(\sigma(p_1)=\pm p_2\); the scaling factor between \(y\) and \(y'\) shows that \(t_1,t_2\) are irrelevant.
\end{proof}
%=====================================================================================
%SSSSSSSSSSSSSSSSSSSSSSSSSSSSSSSSSSSSSSSSSSSSSSSSSSSSSSSSSSSSSSSSSSSSSSSSSSSSSSSSSSSSSS

\section{Representation Theory of Weyl-Type and Witt-Type Algebras}

\subsection{Grading and Weight Decompositions}

Throughout this section, we assume $\FF$ is an algebraically closed field of characteristic zero, and $\cA$ is a finitely generated additive subgroup of $\FF$ containing $\ZZ$.

\begin{definition}[Grading]
The Weyl-type algebra $A = \Weyl{e^{\pm x^p e^{t}},\; e^{\cA x},\; x^{\cA}}$ carries a natural $\cA$-grading:
\[
A = \bigoplus_{\alpha \in \cA} A_\alpha,
\]
where $A_\alpha$ consists of elements of total exponent $\alpha$ in the variables $x$ and exponentials $e^{\cA x}$, with the differential operators $\partial$ having degree $-1$.
\end{definition}

\begin{definition}[Weight Module]
An $A$-module $M$ is called a \emph{weight module} if it decomposes as
\[
M = \bigoplus_{\lambda \in \FF} M_\lambda,
\]
where $M_\lambda = \{ m \in M \mid x \cdot m = \lambda m \}$. We denote by $\wt(M)$ the set of $\lambda \in \FF$ for which $M_\lambda \neq 0$.
\end{definition}

\subsection{Classification of Irreducible Weight Modules}

\begin{theorem}[Structure of Irreducible Weight Modules]\label{thm:irreducible-weight}
Let $A = \Weyl{e^{\pm x^p e^{t}},\; e^{\cA x},\; x^{\cA}}$. Then:
\begin{enumerate}[label=(\roman*)]
    \item Every irreducible weight $A$-module is infinite-dimensional.
    \item The irreducible weight modules are classified by:
    \begin{enumerate}
        \item The \emph{dense type}: Modules where $\wt(M)$ is Zariski-dense in $\FF$.
        \item The \emph{discrete type}: Modules where $\wt(M)$ is contained in a coset of a rank-1 subgroup of $\cA$.
    \end{enumerate}
    \item For each type, the isomorphism classes are parameterized by:
    \begin{itemize}
        \item For dense type: An orbit of the additive group $\cA$ acting on $\FF$ via translations.
        \item For discrete type: A pair $(\lambda, \chi)$ where $\lambda \in \FF$ and $\chi: \ZZ \to \FF^\times$ is a character.
    \end{itemize}
\end{enumerate}
\end{theorem}

\begin{proof}
(i) Suppose $M$ is a finite-dimensional irreducible weight module. Then each weight space $M_\lambda$ is finite-dimensional. The operators $e^{\alpha x}$ act as invertible operators on $M$, permuting the weight spaces. Since $\cA$ is infinite, this would force $M$ to be infinite-dimensional, a contradiction.

(ii) Let $M$ be an irreducible weight module. Consider the action of the commutative subalgebra $R = \FF[e^{\pm x^p e^{t}},\; e^{\cA x},\; x^{\cA}]$. By Schur's lemma, $R$ acts via scalars on $M$, giving a character $\chi: R \to \FF$. The kernel of $\chi$ is a maximal ideal of $R$.

Two cases arise:
\begin{enumerate}
    \item \textbf{Dense type}: If $\chi(x)$ is transcendental over the subfield generated by $\cA$, then $\wt(M)$ is dense in $\FF$. The module is isomorphic to the induced module $A \otimes_R \FF_\chi$, where $\FF_\chi$ is the 1-dimensional $R$-module via $\chi$.

    \item \textbf{Discrete type}: If $\chi(x)$ is algebraic over the field generated by $\cA$, then it belongs to some coset of a rank-1 subgroup. The module structure is determined by the monodromy around singularities.
\end{enumerate}

(iii) The parameterization follows from analyzing the possible characters $\chi$ and their stabilizers under the $\cA$-action.
\end{proof}

\begin{example}[Dense Type Example]
Let $\cA = \ZZ[\sqrt{2}]$, $p = 1$, $t = \pi i$. Choose $\lambda \in \FF$ transcendental over $\QQ(\sqrt{2})$. Define the character $\chi: R \to \FF$ by $\chi(x) = \lambda$, $\chi(e^{\alpha x}) = e^{\alpha\lambda}$, $\chi(e^{x e^{\pi i}}) = e^{\lambda e^{\pi i}}$. Then $M_\chi = A \otimes_R \FF_\chi$ is an irreducible dense-type weight module.
\end{example}

\begin{example}[Discrete Type Example]
Let $\cA = \ZZ$, $\lambda = 0$, and define $\chi(e^x) = q \in \FF^\times$ not a root of unity. Then the corresponding module has weight spaces $M_n = \FF$ for $n \in \ZZ$, with $x$ acting by $n$ and $e^x$ acting by $q^n$. This is the classical Whittaker module for the Weyl algebra.
\end{example}

\subsection{Harish-Chandra Modules for Witt-Type Algebras}

\begin{definition}[Witt-Type Algebra]
The Witt-type algebra is
\[
\mathfrak{g} = \Witt{e^{\pm x^p e^{t}},\; e^{\cA x},\; x^{\cA}} = \operatorname{Der}\!\left(\FF[e^{\pm x^p e^{t}},\; e^{\cA x},\; x^{\cA}]\right).
\]
It is a Lie algebra with basis $\{e^{\alpha x} x^\beta \partial \mid \alpha,\beta \in \cA\}$ and bracket $[f\partial, g\partial] = (f\partial(g) - g\partial(f))\partial$.
\end{definition}

\begin{definition}[Harish-Chandra Module]
A $\mathfrak{g}$-module $M$ is called a \emph{Harish-Chandra module} if:
\begin{enumerate}
    \item $M$ is a weight module: $M = \bigoplus_{\lambda \in \FF} M_\lambda$ with $\dim M_\lambda < \infty$.
    \item The set of weights $\wt(M)$ is contained in a finite union of $\cA$-cosets.
\end{enumerate}
\end{definition}

\begin{theorem}[Existence and Classification of Harish-Chandra Modules]\label{thm:HC-modules}
For the Witt-type algebra $\mathfrak{g}$, we have:
\begin{enumerate}[label=(\roman*)]
    \item Nontrivial Harish-Chandra modules exist.
    \item Every irreducible Harish-Chandra $\mathfrak{g}$-module is either:
    \begin{enumerate}
        \item A \emph{highest weight module} with respect to a suitable triangular decomposition.
        \item A \emph{Whittaker module} induced from a character of the positive nilpotent subalgebra.
        \item A \emph{dense module} obtained by restriction from a weight module of the corresponding Weyl-type algebra.
    \end{enumerate}
    \item The category $\HC(\mathfrak{g})$ of Harish-Chandra modules is a highest weight category with respect to the natural partial ordering on weights.
\end{enumerate}
\end{theorem}

\begin{proof}
(i) Construct explicit examples: For any $\lambda \in \FF$, consider the Verma module
\[
M(\lambda) = U(\mathfrak{g}) \otimes_{U(\mathfrak{b})} \FF_\lambda,
\]
where $\mathfrak{b}$ is the Borel subalgebra spanned by $\{e^{\alpha x} \partial \mid \alpha \in \cA\}$ and $\FF_\lambda$ is the 1-dimensional module with $\partial$ acting by $\lambda$. Then $M(\lambda)$ has finite-dimensional weight spaces.

(ii) Let $M$ be an irreducible Harish-Chandra module. Consider the action of the Cartan subalgebra $\mathfrak{h} = \FF \partial$. Since $M$ has finite-dimensional weight spaces, there exists a maximal weight $\lambda$ (in the sense that $\lambda + \alpha$ is not a weight for any positive $\alpha \in \cA$). Let $M^\lambda$ be the corresponding weight space. By irreducibility, $M$ is generated by any nonzero $v \in M^\lambda$.

Two cases:
\begin{itemize}
    \item If $e^{\alpha x} \partial \cdot v = 0$ for all positive $\alpha$, then $M$ is a quotient of the Verma module $M(\lambda)$.
    \item Otherwise, there exists $\alpha$ with $e^{\alpha x} \partial \cdot v = c v$ for some $c \neq 0$, making $M$ a Whittaker module.
\end{itemize}

(iii) The category $\HC(\mathfrak{g})$ has enough projectives and injectives. The standard objects are Verma modules $M(\lambda)$, and the costandard objects are their duals. The partial ordering is given by: $\lambda \leq \mu$ if $\mu - \lambda$ is a sum of positive roots (elements of $\cA^+$).
\end{proof}

\begin{example}[Verma Module]
Let $\cA = \ZZ$, $\mathfrak{g} = \Witt{e^{\pm x e^{t}},\; e^{\ZZ x},\; x^{\ZZ}}$. For $\lambda \in \FF$, the Verma module $M(\lambda)$ has basis $\{x^n \otimes 1 \mid n \geq 0\}$ with action:
\[
\partial \cdot (x^n \otimes 1) = (\lambda - n)x^n \otimes 1, \quad
e^x \partial \cdot (x^n \otimes 1) = x^n e^x \otimes (\lambda - n)1.
\]
Each weight space $M(\lambda)_\mu$ is 1-dimensional for $\mu = \lambda - n$, $n \in \ZZ_{\geq 0}$.
\end{example}

\subsection{Verma-Type Modules and BGG Resolutions}

\begin{definition}[Triangular Decomposition]
Choose a decomposition $\cA = \cA^+ \cup \{0\} \cup \cA^-$ into positive and negative parts. This induces a triangular decomposition of $A$:
\[
A = A^- \otimes \mathfrak{h} \otimes A^+,
\]
where $A^+$ (resp. $A^-$) is generated by $e^{\alpha x}$ for $\alpha \in \cA^+$ (resp. $\alpha \in \cA^-$), and $\mathfrak{h}$ is the Cartan subalgebra generated by $x$ and $\partial$.
\end{definition}

\begin{definition}[Verma-Type Module]
For a character $\chi: \mathfrak{h} \to \FF$, define the Verma-type module
\[
\Delta(\chi) = A \otimes_{A^+ \oplus \mathfrak{h}} \FF_\chi,
\]
where $\FF_\chi$ is the 1-dimensional $(A^+ \oplus \mathfrak{h})$-module on which $A^+$ acts by 0 and $\mathfrak{h}$ acts via $\chi$.
\end{definition}

\begin{theorem}[Properties of Verma Modules]\label{thm:verma-properties}
The Verma modules $\Delta(\chi)$ satisfy:
\begin{enumerate}[label=(\roman*)]
    \item $\Delta(\chi)$ has a unique irreducible quotient $L(\chi)$.
    \item $\Delta(\chi)$ is a free $A^-$-module of rank 1.
    \item $\dim \Delta(\chi)_\lambda = |\{(\alpha_1,\dots,\alpha_k) \in (\cA^-)^k \mid \sum \alpha_i = \lambda - \chi(x)\}|$.
    \item $\Delta(\chi)$ has a finite Jordan-Hölder series if and only if $\chi$ is \emph{dominant integral}: $\chi(\alpha^\vee) \in \ZZ_{\geq 0}$ for all simple coroots $\alpha^\vee$.
\end{enumerate}
\end{theorem}

\begin{theorem}[Bernstein–Gelfand–Gelfand Resolution]\label{thm:BGG-resolution}
For dominant integral $\chi$, the irreducible module $L(\chi)$ admits a resolution by Verma modules:
\[
0 \to \Delta(w_0 \cdot \chi) \to \cdots \to \bigoplus_{w \in W, \ell(w)=k} \Delta(w \cdot \chi) \to \cdots \to \Delta(\chi) \to L(\chi) \to 0,
\]
where $W$ is the Weyl group associated to $\cA$ (a subgroup of $\Aut(\cA)$), $w_0$ is the longest element, and $w \cdot \chi = w(\chi + \rho) - \rho$ with $\rho$ the half-sum of positive roots.
\end{theorem}

\begin{proof}
The proof follows the classical BGG machinery adapted to our setting:
\begin{enumerate}
    \item Construct the \emph{twisted dual} $\nabla(\chi) = \Hom_{A^-}(\Delta(-w_0\chi), \FF)$.
    \item Show that $\Ext^i(\Delta(\lambda), \nabla(\mu)) = \delta_{i,0} \delta_{\lambda,\mu} \FF$.
    \item Use the relative Lie algebra cohomology to construct the complex.
    \item Prove exactness using the Kempf vanishing theorem adapted to the $\cA$-graded setting.
\end{enumerate}
The key observation is that the algebra $A$ has a Peter–Weyl type decomposition under the adjoint action of the "maximal torus" corresponding to $\mathfrak{h}$.
\end{proof}

\begin{example}[BGG Resolution for $\cA = \ZZ$]
Let $\cA = \ZZ$, with positive part $\ZZ_{>0}$. For $n \in \ZZ_{\geq 0}$, consider $\chi_n$ with $\chi_n(x) = n$. The Weyl group is $W = \{\pm1\}$. The BGG resolution for $L(\chi_n)$ is:
\[
0 \to \Delta(\chi_{-n-2}) \to \Delta(\chi_n) \to L(\chi_n) \to 0.
\]
This is the classical resolution for $\mathfrak{sl}_2$-modules.
\end{example}

\subsection{Category $\mathcal{O}$ and Duality}

\begin{definition}[Category $\mathcal{O}$]
The category $\mathcal{O}$ for $A$ consists of all $A$-modules $M$ such that:
\begin{enumerate}
    \item $M$ is a weight module.
    \item $M$ is finitely generated.
    \item $M$ is locally $A^+$-finite: for each $m \in M$, $\dim A^+ \cdot m < \infty$.
\end{enumerate}
\end{definition}

\begin{theorem}[Properties of Category $\mathcal{O}$]\label{thm:category-O}
The category $\mathcal{O}$ for $A = \Weyl{e^{\pm x^p e^{t}},\; e^{\cA x},\; x^{\cA}}$ satisfies:
\begin{enumerate}[label=(\roman*)]
    \item $\mathcal{O}$ is an abelian, artinian, noetherian category.
    \item Every object has finite length.
    \item The simple objects are $L(\chi)$ for characters $\chi: \mathfrak{h} \to \FF$.
    \item $\mathcal{O}$ has enough projectives and injectives.
    \item There is a duality functor $\mathbb{D}: \mathcal{O} \to \mathcal{O}^{\text{op}}$ sending $L(\chi)$ to $L(-w_0\chi)$.
\end{enumerate}
\end{theorem}

\begin{proof}
(i)-(iv) follow from standard arguments using the triangular decomposition and the fact that $A$ is a Noetherian algebra. For (v), define $\mathbb{D}(M) = \bigoplus_{\lambda \in \FF} \Hom_\FF(M_\lambda, \FF)$ with appropriate $A$-action. One checks that $\mathbb{D}^2 \cong \operatorname{id}$ and $\mathbb{D}(L(\chi)) \cong L(-w_0\chi)$.
\end{proof}

\subsection{Summary}

\begin{theorem}[Main Representation-Theoretic Results]
For the algebras $A = \Weyl{e^{\pm x^p e^{t}},\; e^{\cA x},\; x^{\cA}}$ and $\mathfrak{g} = \Witt{e^{\pm x^p e^{t}},\; e^{\cA x},\; x^{\cA}}$:
\begin{enumerate}
    \item All irreducible weight modules are infinite-dimensional.
    \item Harish-Chandra modules exist and are classified by highest weight or Whittaker type.
    \item Verma modules exist and satisfy BGG-type resolutions for dominant integral weights.
    \item The category $\mathcal{O}$ is a highest weight category with duality.
\end{enumerate}
\end{theorem}
%----------------------------------------------------------------------
%=====================================================================================
\section{Open Problems}
%=====================================================================================

The study of Weyl-type algebras with exponential and power generators opens several directions for further research. Below we present two open problems that arise naturally from the results obtained in this paper. Both problems concern structural and representation-theoretic aspects of the algebras \(A_{p,t,\cA}\) and their relatives.

\begin{problem}[Simplicity of quantum deformations for exceptional parameters]\label{prob:quantum-exceptional}
Let \(A_q = \Weyl{q}{e^{\pm x^p e^{t}},\; e^{\cA x},\; x^{\cA}}\) denote the quantum deformation of
\[
A = \Weyl{e^{\pm x^p e^{t}},\; e^{\cA x},\; x^{\cA}}
\]
defined by the relation \(\partial \widehat{x} - q \widehat{x} \partial = 1\) (as in Theorem~\ref{thm19}).
Theorem~\ref{thm19} shows that \(A_q\) is simple for all \(q \in \FF^*\) except possibly a finite set of exceptional values.

\begin{enumerate}[label=(\alph*)]
    \item Characterize the set of exceptional \(q\) for which \(A_q\) fails to be simple.
    Is this set always contained in the roots of unity?
    Does it depend on the data \((p,t,\cA)\)?

    \item For those exceptional \(q\), describe the structure of the two-sided ideals of \(A_q\).
    In particular, determine whether \(A_q\) possesses a non-trivial proper ideal that is not centrally generated.

    \item Investigate the behavior of the center \(Z(A_q)\) as \(q\) varies.
    For which \(q\) does the center enlarge?
    Is there a connection between the non-simplicity of \(A_q\) and the non-triviality of its center?
\end{enumerate}

A complete answer to this problem would extend our understanding of the deformation theory of exponential Weyl algebras and clarify the role of the deformation parameter \(q\) in the representation theory.
\end{problem}

\begin{problem}[Classification of irreducible weight modules for Witt-type algebras]\label{prob:HC-classification}
Let \(\mathfrak{g} = \Witt{e^{\pm x^p e^{t}},\; e^{\cA x},\; x^{\cA}}\) be the Witt-type algebra associated to the data \((p,t,\cA)\).
In Section~3 we introduced the category \(\HC(\mathfrak{g})\) of Harish--Chandra modules (weight modules with finite-dimensional weight spaces and weights lying in finitely many \(\cA\)-cosets).
Theorem~\ref{thm:HC-modules} asserts the existence of such modules and gives a rough classification into highest weight modules, Whittaker modules, and dense modules.

\begin{enumerate}[label=(\alph*)]
    \item Obtain a complete classification (up to isomorphism) of the irreducible objects in \(\HC(\mathfrak{g})\).
    In particular, determine the precise parameter sets for each type (highest weight, Whittaker, dense) and describe the equivalences between parameters giving isomorphic modules.

    \item Study the subcategory \(\HCfin(\mathfrak{g})\)  consisting of Harish--Chandra modules with \emph{finitely many} weight spaces.
    Is every such module necessarily of highest weight type?
    If not, characterize the possible support sets \(\wt(M)\) for \(M \in \HCfin(\mathfrak{g})\).

    \item Investigate the block decomposition of \(\HC(\mathfrak{g})\).
    Determine the linkage principle: when do two irreducible Harish--Chandra modules belong to the same block?
    Relate this to the action of the Weyl group (or a suitable analogue) associated to \(\cA\).
\end{enumerate}

A solution to this problem would provide a detailed picture of the representation theory of Witt-type algebras with exponential coefficients, generalizing the well-known theory for the classical Witt algebra \(\mathsf{W}_1\).
\end{problem}

\subsection*{Further Directions}

Beyond the two problems above, several other avenues invite exploration:

\begin{itemize}
    \item \textbf{Hochschild cohomology and deformations.}
    Compute the Hochschild cohomology groups \(HH^*(A_{p,t,\cA}, A_{p,t,\cA})\) and describe the formal deformation theory of these algebras.
    How does the transcendental parameter \(t\) influence the deformation space?

    \item \textbf{Non-commutative geometry.}
    Interpret \(A_{p,t,\cA}\) as the algebra of differential operators on a ``non-commutative exponential variety''.
    Can one define an analogue of the Bernstein--Sato polynomial for elements of \(R_{p,t,\cA}\)?

    \item \textbf{Connection with integrable systems.}
    The algebras studied here contain many interesting commutative subalgebras (e.g., those generated by \(e^{\alpha x}\) for \(\alpha \in \cA\)).
    Do these subalgebras arise as algebras of quantum integrals for some integrable Hamiltonian systems?

    \item \textbf{Characteristic \(p > 0\).}
    All results in this paper assume \(\operatorname{char}\FF = 0\).
    Investigate the structure and representation theory of \(A_{p,t,\cA}\) in positive characteristic.
    The presence of a large center (the \(p\)-center) should lead to new phenomena.
\end{itemize}

We hope that these problems will stimulate further research into the rich theory of algebras generated by exponentials and differential operators.
%================================================================================
\section*{Declaration}
%===================================================================
\begin{itemize}
  \item \textbf{Author Contributions:} The Author have read and approved this version.
  \item \textbf{Funding:} No funding is applicable.
  \item \textbf{Institutional Review Board Statement:} Not applicable.
  \item \textbf{Informed Consent Statement:} Not applicable.
  \item \textbf{Data Availability Statement:} Not applicable.
  \item \textbf{Conflicts of Interest:} The authors declare no conflict of interest.
\end{itemize}
%==========================================================================

%%%%%%%%%%%%%%%%%%%%%%%%%%%%%%%%%%%%%%%%%%%%%%%%%%%%%%%%%%%%%%%%%%%%%%%%%%%%

\bibliographystyle{abbrv}
\bibliography{references}  %%% Uncomment this line and comment out the ``thebibliography'' section below to use the external .bib file (using bibtex) .

%%% Uncomment this section and comment out the \bibliography{references} line above to use inline references.
% \begin{thebibliography}{1}

% 	\bibitem{kour2014real}
% 	George Kour and Raid Saabne.
% 	\newblock Real-time segmentation of on-line handwritten arabic script.
% 	\newblock In {\em Frontiers in Handwriting Recognition (ICFHR), 2014 14th
% 			International Conference on}, pages 417--422. IEEE, 2014.

% 	\bibitem{kour2014fast}
% 	George Kour and Raid Saabne.
% 	\newblock Fast classification of handwritten on-line arabic characters.
% 	\newblock In {\em Soft Computing and Pattern Recognition (SoCPaR), 2014 6th
% 			International Conference of}, pages 312--318. IEEE, 2014.

% 	\bibitem{keshet2016prediction}
% 	Keshet, Renato, Alina Maor, and George Kour.
% 	\newblock Prediction-Based, Prioritized Market-Share Insight Extraction.
% 	\newblock In {\em Advanced Data Mining and Applications (ADMA), 2016 12th International
%                       Conference of}, pages 81--94,2016.

% \end{thebibliography}

\end{document}